\documentclass[hypertex,11pt]{article}
\usepackage{amsfonts,amsmath,amscd,amssymb,amsthm}
\usepackage{fullpage}
\usepackage[backrefs, alphabetic,initials]{amsrefs}
\usepackage[all]{xy}
\usepackage{hyperref} 
\usepackage{mathrsfs}

\def\convf{\hbox{\space \raise-2mm\hbox{$\textstyle      \bigotimes \atop \scriptstyle \omega$} \space}}
\def\0{{\bar 0}}
\def\1{{\bar 1}}

\def\Z{{\mathbb Z}}

\def\N{{\mathbb N}}

\def\V{\mathscr{V}}

\def\ad{{\operatorname{ad \;}}}

\def\trp{{\operatorname{transpose}}}

\def\hgt{{\operatorname{ht}}}

\def\gra{{\operatorname{gr}\;}}
\def\gra{{\operatorname{gr}\;}}

\def\ch{{\operatorname{ch}\:}}

\def\Hom {{\operatorname{Hom}}}

\def\ann {{\operatorname{ann}}}

\newcommand{\ttk}{\mathtt{k}}

\newcommand{\gL}{\Lambda}

\newcommand{\itema}{\item[{{\rm(a)}}]}

\newcommand{\itemb}{\item[{{\rm(b)}}]}
\newcommand{\itemc}{\item[{{\rm(c)}}]}

\newcommand{\da}{\downarrow}

\newcommand{\noi}{\noindent}
\newcommand{\ga}{\alpha}
\newcommand{\gb}{\beta}
\newcommand{\gc}{\gamma}

\newcommand{\Gd}{\Delta}
\newcommand{\gd}{\delta}

\newcommand{\gs}{\sigma}

\newcommand{\cg}{{\rm Cdeg}}
\newcommand{\gt}{\tau}
\newcommand{\gz}{\zeta}
\newcommand{\gl}{\lambda}

\newcommand{\gr}{\rho}

\newcommand{\gep}{\epsilon}
\newcommand{\gth}{\theta}

\newcommand{\gve}{\varepsilon}

\newcommand{\fg}{\mathfrak{g}}\newcommand{\fgl}{\mathfrak{gl}}\newcommand{\so}{\mathfrak{so}}
\newcommand{\fsl}{\mathfrak{sl}}\newcommand{\fpsl}{\mathfrak{psl}}\newcommand{\osp}{\mathfrak{osp}}

\newcommand{\fr}{\mathfrak{r}}

\newcommand{\fh}{\mathfrak{h}}

\newcommand{\fb}{\mathfrak{b}}

\newcommand{\fn}{\mathfrak{n}}

\newcommand{\fsp}{\mathfrak{sp}}
\newcommand{\fk}{\mathfrak{k}}
\newcommand{\fp}{\mathfrak{p}}

\newcommand{\sfb}{\small{\mathfrak{b}}}

\newcommand{\sfg}{\small{\mathfrak{g}}}

\newcommand{\ff}{\footnote}
\newfont{\eufm}{eufm10 scaled\magstep1}

\newcommand{\cO}{\mathcal{O}}

\newcommand{\cH}{\mathcal{H}}

\newcommand{\bbZ}{\mathbb{Z}}

\newcommand{\ey}{\end{eqnarray}}
\newcommand{\by}{\begin{eqnarray}}

\newcommand{\bco}{\begin{conjecture}}
\newcommand{\ba}{\begin{alg}}
\newcommand{\ea}{\end{alg}}
\newcommand{\eco}{\end{conjecture}}
\newcommand{\bpf}{\begin{proof}}
\newcommand{\epf}{\end{proof}}
\newcommand{\bt}{\begin{theorem}}

\newcommand{\et}{\end{theorem}}

\newcommand{\br}{\begin{rem}}
\newcommand{\er}{\end{rem}}
\newcommand{\brs}{\begin{rems}}
\newcommand{\ers}{\end{rems}}
\newcommand{\bi}{\begin{itemize}}
\newcommand{\bl}{\begin{lemma}}
\newcommand{\bsul}{\begin{sublemma}}
\newcommand{\esul}{\end{sublemma}}
\newcommand{\bp}{\begin{proposition}}
\newcommand{\be}{\begin{equation}}
\newcommand{\bc}{\begin{corollary}}
\newcommand{\bexs}{\begin{examples}}
\newcommand{\eexs}{\end{examples}}
\newcommand{\bexa}{\begin{example}}
\newcommand{\eexa}{\end{example}}
\newcommand{\bex}{\begin{exercise}}
\newcommand{\eex}{\end{exercise}}
\newcommand{\btab}{\begin{tab}}
\newcommand{\etab}{\end{tab}}
\newcommand{\ei}{\end{itemize}}
\newcommand{\el}{\end{lemma}}
\newcommand{\ep}{\end{proposition}}
\newcommand{\ee}{\end{equation}}
\newcommand{\ec}{\end{corollary}}
\newcommand{\Bc}{\begin{center}}
\newcommand{\Ec}{\end{center}}
\newcommand{\bh}{\begin{hyp}}
\newcommand{\eh}{\end{hyp}}
\newcommand{\bhs}{\begin{hyps}}
\newcommand{\ehs}{\end{hyps}}
\newcommand{\bd}{\begin{dfn}}
\newcommand{\ed}{\end{dfn}}

\begin{document}
\title{Table of Contents}

\newtheorem{thm}{Theorem}[section]
\newtheorem{hyp}[thm]{Hypothesis}
 \newtheorem{hyps}[thm]{Hypotheses}
  \newtheorem{rems}[thm]{Remarks}

\newtheorem{conjecture}[thm]{Conjecture}
\newtheorem{theorem}[thm]{Theorem}
\newtheorem{theorem a}[thm]{Theorem A}
\newtheorem{example}[thm]{Example}
\newtheorem{examples}[thm]{Examples}
\newtheorem{corollary}[thm]{Corollary}
\newtheorem{rem}[thm]{Remark}
\newtheorem{lemma}[thm]{Lemma}
\newtheorem{cor}[thm]{Corollary}
\newtheorem{proposition}[thm]{Proposition}
\newtheorem{exs}[thm]{Examples}
\newtheorem{ex}[thm]{Example}
\newtheorem{exercise}[thm]{Exercise}
\numberwithin{equation}{section}%
\setcounter{part}{0}
\newcommand{\rra}{\longleftarrow}
\newcommand{\lra}{\longrightarrow}
\newcommand{\dra}{\Rightarrow}
\newcommand{\dla}{\Leftarrow}
\newtheorem{Thm}{Main Theorem}


\newtheorem*{thm*}{Theorem}
\newtheorem{lem}[thm]{Lemma}
\newtheorem*{lem*}{Lemma}
\newtheorem*{prop*}{Proposition}
\newtheorem*{cor*}{Corollary}
\newtheorem{dfn}[thm]{Definition}
\newtheorem*{defn*}{Definition}
\newtheorem{notadefn}[thm]{Notation and Definition}
\newtheorem*{notadefn*}{Notation and Definition}
\newtheorem{nota}[thm]{Notation}
\newtheorem*{nota*}{Notation}
\newtheorem{note}[thm]{Remark}
\newtheorem*{note*}{Remark}
\newtheorem*{notes*}{Remarks}
\newtheorem{hypo}[thm]{Hypothesis}
\newtheorem*{ex*}{Example}
\newtheorem{prob}[thm]{Problems}
\newtheorem{conj}[thm]{Conjecture}

\title{ Coefficients of \v Sapovalov elements for simple Lie algebras and contragredient Lie superalgebras.}
\author{Ian M. Musson\ff{Research partly supported by  NSA Grant H98230-12-1-0249.} \\Department of Mathematical Sciences\\
University of Wisconsin-Milwaukee\\ email: {\tt
musson@uwm.edu}}
\maketitle
\begin{abstract}
We provide upper bounds on the degrees of the coefficients of \v{S}apovalov elements for a simple Lie algebra.  If $\fg$ is a contragredient Lie superalgebra
and $\gc$ is a positive isotropic root of $\fg,$ we prove the existence and uniqueness of the  \v Sapovalov element for $\gc$ and we obtain upper bounds on the degrees of their coefficients. For type A Lie superalgebras we give a closed formula for \v Sapovalov elements.
 We also explore the
 behavior of \v Sapovalov elements
  when the Borel subalgebra is changed, and the survival of \v{S}apovalov elements in factor modules of Verma modules.
\end{abstract}
\Bc Table of Contents. \Ec
\ref{uv.1}. {Introduction.} \\
\ref{s.1}. {Uniqueness of \v Sapovalov elements.} \\
\ref{s.4}. {Outline of the Proof and Preliminary Lemmas.} \\
\ref{s.5}. {Proof of Theorem \ref{Shap}.}  \\
\ref{s.51}. {Proof of Theorem \ref{aShap}.}\\
\ref{zzprod}. {Powers of  \v Sapovalov elements.}\\
  \ref{cosp}. {An (ortho) symplectic example.} \\
\ref{s.8}. {The Type A Case.} \\
\ref{surv}. {{Survival of \v{S}}apovalov elements in factor modules.}\\ 
\ref{sscbs}. {Changing the Borel subalgebra.} \\
\section{Introduction.} \label{uv.1}
Throughout this paper we work over an  algebraically closed field $\ttk $
  of characteristic zero. If  $\fg$ is a simple Lie algebra
necessary and sufficient conditions for the existence of a non-zero homomorphism from $M(\mu)$ to $M(\gl)$ can be obtained by combining work of Verma
\cite{Ve} with work of Bernstein, Gelfand and  Gelfand \cite{BGG1}, \cite{BGG2}.  Such maps can be described explicitly in terms of certain elements introduced by N.N.\v Sapovalov in \cite{Sh}.
Verma modules are fundamental objects in the study of category $\cO$, a study that has blossomed into
an extremely rich theory in the years since these early papers appeared.  Highlights include the Kazhdan-Lusztig conjecture \cite{KaLu},  \cite{BB}, \cite{BK }, the work of
Beilinson, Ginzburg, and Soergel on Koszul duality \cite{BGS}, and more recent results on categorification see   \cite{H2} and \cite{Ma} for more details.
A refined version of the Kazhdan-Lusztig conjecture giving the composition multiplicities in successive quotients of the Jantzen filtration conjecture was 
proven in \cite{BBJ} by showing that the localization functor sends the  Jantzen filtration to the weight
 filtration on perverse sheaves.  This also established a conjecture of  Jantzen on a compatibility  property of the filtration from \cite{J1} with Verma submodules.
\\ \\
Recently significant advances have been made in the study of the category $\cO$ for classical simple Lie superalgebras using a variety of techniques.  After the early work of Kac  \cite{K}, \cite{Kac3} the first major advance was made by Serganova who used geometric techniques to obtain a character formula for Kac modules over $\fgl(m,n)$, \cite{S2}.  The next development was Brundan's approach to the same problem using a combination of algebraic and  combinatorial  techniques.  In his seminal paper  \cite{Br} also introduced a Fock space representation $ \mathcal{T}^{m|n} :=
{\bigotimes}^m \V^* \otimes {\bigotimes}^n \V$  of the quantized enveloping algebra $\mathcal{U}=
U_q(\fgl_\infty)$
where  $\V$ is the natural representation of $\mathcal{U},$ and $\V^*$ is its restricted dual.
He then introduced monomial
and  canonical bases for $\mathcal{T}^{m|n}$ , and
using the transition matrices between these matrices defined polynomials which have
become known as Brundan-Kazhdan-Lusztig polynomials.
\\ \\
  Brundan then made the extraordinary conjecture that the values at $q = 1$  of these polynomials solve the multiplicity problem of for composition factors of Verma modules.  For $\fgl(m)$ this is equivalent to the Kazhdan-Lusztig conjecture. Brundan's conjecture was later confirmed
 by Cheng, Lam and Wang \cite{CLW},  exploiting connections with
 super-duality. Super-duality connects the parabolic category $\cO$ for $\fgl(m,n)$ to a corresponding parabolic category for $\fgl(m+n)$. The authors later extended this connection to the orthosymplectic case \cite{CLW2}, \cite{CW}.
 A new proof of the conjecture was provided by
  Brundan, Losev and Webster \cite{BLW}, at the same time showing that any integral block
   of the category $\cO$ for $\fgl(m,n)$ has a graded lift which is Koszul, see the recent survey article \cite{B} by Brundan for these developments.
 \\ \\The \v Sapovalov determinant, also introduced in \cite{Sh} has been developed in a variety of contexts, such as Kac-Moody algebras \cite{KK}, quantum groups \cite{Jo1},
and Lie superalgebras \cite{G4}, \cite{G}, \cite{G2}.  However neither \v Sapovalov elements
nor the  Jantzen filtration have received much attention
for classical simple Lie superalgebras.
The purpose of this paper is to initiate the study of
\v Sapovalov elements  in the super case. New phenomena arise due to the presence of isotropic roots.   A sequel will focus on the Jantzen filtration  and sum formula \cite{M17}.
\\ \\
 Let $\fg= \fg(A,\gt)$ be a finite dimensional contragredient Lie superalgebra with Cartan subalgebra $\fh$, and set of simple roots $\Pi$.
 The superalgebras $\fg(A,\gt)$ coincide with the basic classical simple Lie superalgebras, except that instead of $\fpsl(n,n)$ we obtain $\fgl(n,n).$
 Implicit in the definition of the $\fg(A,\gt)$ is a preferred Borel subalgebra.
 Let $\Delta^{+}$ and  %
\be \label{gtri}\mathfrak{g} = \mathfrak{n}^- \oplus \mathfrak{h}
\oplus \mathfrak{n}^+\ee
be the set of positive roots  containing $\Pi$, and  the corresponding triangular decomposition  of $\fg$ respectively. We use the Borel subalgebras $\mathfrak{b} =  \mathfrak{h}
\oplus \mathfrak{n}^+$ and
$\fb^- =
\mathfrak{n}^- \oplus \mathfrak{h}$.
 The Verma module
$M(\gl)$ with highest weight $\gl \in \fh^*$, and highest weight vector
$v_\lambda$
is induced from $\mathfrak{b}$. Suppose that $\gc$  is a positive root, and $m$ is a  positive   integer.  The \v Sapovalov element
$\theta_{\gamma,m}$
 corresponding to the pair $(\gc, m)$ has the form \be \label{rat}
\theta_{\gamma,m} = \sum_{\pi \in {\overline{\bf P}}(m\gamma)} e_{-\pi}
H_{\pi},\ee where $H_{\pi} \in U({\mathfrak h})$,
and has the property that if $\gl$ lies on a certain hyperplane then
$\theta_{\gc, m}v_{\gl}$ is a highest weight vector  in $M(\gl),$ see (\ref{boo}).
In (\ref{rat}) the sum is indexed by the set
${\overline{\bf P}}(m\gamma)$ of partitions of $m\gc  $, the
$e_{-\pi}$ with $\pi \in {\overline{\bf P}}(m\gamma)$ form a basis for the weight space $U({\mathfrak \fn^-})^{-m\gc}$, and the coefficients
$H_{\pi}$ are in $U({\mathfrak h})$.
We normalize  $\theta_{\gamma,m}$ so that for a certain ${\pi^0}
\in {\overline{\bf P}}(m\gamma)$, the coefficient $H_{\pi^0}$ is equal to 1.
This guarantees that $\theta_{\gc, m}v_\gl$ is never zero.
\\ \\
  The main results in this paper 
 give bounds on the degrees of the coefficients $H_{\pi}$ in (\ref{rat}).  There is always a unique coefficient  of highest degree, and we determine the leading term of this coefficient up to a scalar multiple. These results appear to be new even for simple Lie algebras.
The exact form of the coefficients depends on the way the positive roots are ordered.  Nevertheless they seem to have interesting properties both combinatorially and from the point of view of representation theory.  For example with a suitable ordering the coefficients are often products of linear factors and the vanishing of  these factors has an interpretation in terms of representation theory.
\\ \\
 The existence of a unique coefficient  of highest degree is useful in the construction of  some new highest weight modules $M_\gc(\gl)$, where $\gc$ is an isotropic root and $(\gl+\gr,\gc) =0$, \cite{M17}.  This results in an improvement in the Jantzen sum formula from \cite{M} Theorem 10.3.1.  The module  $M_\gc(\gl)$ has character  $\epsilon^{\lambda}p_\gc$ see (\ref{jfn}) for notation, so this leads to a formula where both sides are sums of characters in the category $\cO.$
 \\ \\
 \v Sapovalov elements corresponding to non-isotropic roots for  a basic classical simple  Lie superalgebra
were constructed  in \cite{M} Chapter 9.  This closely parallels the semisimple
case.
Properties of the coefficients of these elements were announced in \cite{M} Theorem 9.2.10. However the bounds on the degrees of the coefficients claimed in \cite{M} are incorrect
if $\Pi$ contains a non-isotropic odd root.
They are corrected by Theorem \ref{aShap}.\\ \\
I would like to thank Jon Brundan for suggesting the use of noncommutative determinants to write \v Sapovalov elements in Section \ref{s.8}, and raising the possibility of using Theorem \ref{shgl} to prove Theorem \ref{shapel}. I also thank Kevin Coulembier for some helpful conversations.
\subsection{\bf Preliminaries.}  \label{sss7.1}
\noi
We use the definition of partitions from \cite{M} Remark 8.4.3. Set $Q^+=\sum_{\ga\in \Pi} \N\ga$.
If $\eta \in Q^+$, a {\it
partition} of $\eta$ is a map
$\pi: \Delta^+ \longrightarrow
\mathbb{N} $ such that
$\pi(\alpha) = 0$ or $1$ for all isotropic roots $\alpha$,
 $\pi(\alpha) = 0$ for all even roots $\alpha$ such that $\ga/2$ is a root, and
\[ \sum_{\alpha \in {{\Delta^+}}} \pi(\alpha)\alpha = \eta.\]
For $\eta \in Q^+$ , we denote by $\bf{\overline{P}(\eta)}$ the set of partitions of $\eta$. If $\pi \in \bf{\overline{P}(\eta)}$  the {\it degree} of $\pi$ is defined to be $|\pi| = \sum_{\alpha \in \Delta^+} \pi(\alpha).$\\ \\
\noi Fix a non-degenerate invariant symmetric bilinear form $(\;,\;)$ on $\fh^*$, and for all $\ga \in \fh^*$, let $h_\ga \in \fh$ be the unique element such that $(\ga,\gb) = \gb(h_\ga)$ for all $\gb \in \fh^*$. Then for all $\alpha
\in \Delta^+$, choose elements $e_{\pm \alpha} \in
\mathfrak{g}^{\pm \alpha}$
 such that
\[ [e_{\alpha}, e_{-\alpha}] = h_{\alpha}.\]
\noi  Fix an ordering on the set $\Delta^+$, and for $\pi$ a partition,
set
\[ e_{-\pi} = \prod_{\alpha \in \Delta^+} e^{\pi (\alpha)}_{-\alpha},\]
the product being taken with respect to this order.
Then the elements $e_{- \pi},$ with
$\pi \in \bf{\overline{P}}(\eta)$ form a basis of
$U(\mathfrak{n}^-)^{- \eta}.$  \noi For a non-isotropic root $\ga,$ we set $\alpha^\vee = 2\alpha /
(\alpha, \alpha)$, and denote the reflection corresponding to $\alpha$ by $
s_\ga.$
As usual the  Weyl group  $W$ is the subgroup of $GL({\fh}^*)$
generated by all such reflections.  For $u \in W$ set
 $$ N(u) = \{ \alpha \in \Delta_0^+ | u \alpha < 0 \},\qquad \ell(u) = |N(u)|.$$
We use the following well-known fact several times.
\bl
If $w  = s_\ga u$ with  $\ell(w)>\ell(u)$ and $\ga$ is a simple non-isotropic root, then we have a disjoint union
\begin{eqnarray} \label{Nw}
N(w^{-1}) = s_\ga N(u^{-1}) \cup \{\ga\}.
\end{eqnarray}
\el\bpf See for example \cite{H3} Chapter 1.\epf
\noi Set $$\gr_0({\fb})=\frac{1}{2}\sum_{\ga \in \Gd_0^+}\ga , \quad
\gr_1({\fb})=\frac{1}{2}\sum_{\ga \in \Gd_1^+}\ga ,\quad  \gr({\fb})=
\gr_0({\fb}) -\gr({\fb_1}).$$
Except in Section \ref{sscbs} we work with a fixed Borel subalgebra, and if this is the case we set $\gr_i =
\gr_i({\fb})$ for $i=1,2$ and $\gr=\gr({\fb}).$

\subsection{Main Results.} \label{s.2}
Fix a positive root $\gc$  and a  positive   integer $m$. 
Let
$\pi^0 \in {\overline{\bf P}}(m\gc)$ be the unique partition of $m\gc$ such
that $\pi^0(\ga) = 0$ if $\ga \in \Delta^+ \backslash \Pi.$ The partition $m\pi^{\gamma}$ of $m\gc$ is given by $m\pi^{\gamma}(\gc)=m,$ and $m\pi^{\gamma}(\ga)= 0$ for all positive roots $\ga$ different from $\gc.$
We say that $\gth \in U({\mathfrak b}^{-} )^{- m\gc}$ is a 
{\it \v Sapovalov element for the pair} $(\gc,m)$ if it has the form (\ref{rat}) 
with $H_{\pi^0} = 1,$  and
\be \label{boo} e_{\ga} \theta \in
U({\mathfrak g})(h_{\gc} + \rho(h_{\gc})
-m(\gc,\gc)/2)+U({\mathfrak g}){\mathfrak n}^+ , \; \rm{
for \; all }\;
\ga \in \Delta^+.
 \ee
For a semisimple Lie algebra, the existence of such elements was shown by \v Sapovalov, \cite{Sh} Lemma 1.
Let \[{\mathcal H}_{\gc, m} = \{ \lambda \in  {\mathfrak h}^*|(\lambda + \rho, \gc) = m(\gc, \gc)/2  \},\] and let ${\mathcal I}(\mathcal H_{\gc, m})$ be the ideal of $S(\fh)$ consisting of  functions vanishing on $\mathcal{H}_{\gamma,m}$. Thus the ring of regular functions on ${\mathcal H}_{\gc, m}$ is
\[
 S(\fh)/(h_\gamma + \rho(h_\gamma) - m(\gamma, \gamma)/2).\]
\noi
Note that if $\gl \in \mathcal H_{\gc, m}$, and $\gth$ satisfies (\ref{boo}) then $\theta v_\lambda$ is a highest weight vector of weight $\lambda -m\gc$ in
$M(\gl)$. Now for $\mu, \gl \in \fh^*$ we have
\be \label{kit}\dim \Hom_\fg(M(\mu),M(\gl))\le 1,\ee
by \cite{D} Theorem 7.6.6,
and it follows that the \v Sapovalov element $\theta = \theta_{\gc, m}$ for the pair $(\gc,m)$ is unique modulo the left ideal
$U({\mathfrak b}^{-} ){\mathcal I}(\mathcal H_{\gc, m})$.
\\ \\ \noi
A finite dimensional contragredient Lie superalgebra $\fg$ has, in general several conjugacy classes of Borel subalgebras, and this both complicates and enriches the representation theory of $\fg$.  The complications are partially resolved by at first fixing a Borel subalgebra (or equivalently a basis of simple roots for $\fg$) with special properties.  Later we study the effect of changing the Borel subalgebra.
The first explicit description of the system of roots and possible Dynkin-Kac diagrams, up to conjugacy by the Weyl group, was given by Kac in \cite{K} sections 2.5.4 and 2.5.5.
We remark that there are some omissions on these lists.  The corrected lists appear in \cite{FSS} (and elsewhere).
Note that if $\fg = \osp(2m,2n)$ there are some diagrams that correspond to two sets of simple roots.  The corresponding Borel subalgebras are conjugate under an outer (diagram) automorphism of $\fg$, \cite{M} Corollary 5.5.13.\\ \\
In \cite{K} Table VI Kac gave a particular diagram in each case that we will call {\it distinguished.} The corresponding set of simple roots and Borel subalgebra are also called distinguished. The distinguished Borel subalgebra contains at most one simple isotropic root vector.  Unless $\fg = \osp(1,2n)$, $\fg = \osp(2,2n)$ there is exactly one other Borel subalgebra with this property up to conjugacy in Aut $\fg$. A representative of this class (and its set of simple roots) will be called {\it anti-distinguished.} If $\fg = \osp(1,2n)$  there is only one conjugacy class, while if $\fg = \osp(2,2n)$ there is a basis containing exactly two isotropic roots $\gd_n\pm \gep_1$ (in the notation of Kac). In the latter case we call this basis
anti-distinguished.
  Suppose $\Pi_{\rm nonisotropic},$ (resp. $\Pi_{\rm even}$) be the set of nonisotropic (resp. even) simple roots, and let $W_{\rm nonisotropic}$ (resp. $W_{\rm even}$) be
the subgroup of $W$ generated by the reflections $s_\ga,$ where $\ga
\in \Pi_{\rm nonisotropic}$ (resp. $\ga \in \Pi_{\rm even}$).  Consider the following hypotheses.
\be \label{i} \mbox{ The set of simple roots of } \Pi \mbox{ is  either distinguished or anti-distinguished.} \ee
\be \label{iii} \gc = w\gb \mbox{ for a simple root } \gb \mbox{ and } w \in W_{\rm even}.\ee
\be \label{ii}
\gc = w\gb \mbox{ for a simple root } \gb \mbox{ and } w \in W_{\rm nonisotropic}.\ee
When (\ref{i}) and either of (\ref{iii}) or (\ref{ii}) holds we always assume that $\ell(w)$ is minimal, and for $\alpha \in N(w^{-1}),$ we define $q(w,\ga) = (w\gb, \ga^\vee).$ If $\Pi = \{\ga_i|i   = 1, \ldots, t \}$ is the set of simple roots, and $\gc = \sum_{i=1}^t a_i\ga_i,$ then the {\it height} $\hgt \gc$ of $\gc$ is defined to be
$\hgt \gc = \sum_{i=1}^t a_i$.
\begin{theorem} \label{Shap}
Suppose $\fg$ is semisimple or a contragredient Lie superalgebra, and $\gc$ is a positive root such that $(\ref{i})$ and $(\ref{iii})$ hold.  If $\gc$ is isotropic assume that $m=1.$ Then
there exists a \v Sapovalov  element
$\theta_{\gamma,m} \in U({\mathfrak b}^{-} )^{- m\gamma}$, which is unique modulo the left ideal
$U({\mathfrak b}^{-} ){\mathcal I}(\mathcal H_{\gc})$, and
the coefficients of $\theta_{\gc, m}$ satisfy
\be \label{x1} |\pi|+\deg  H_{\pi} \le m\hgt \gc, \ee and
\be\label{hig} H_{m\pi^{\gamma}} \mbox{
has leading term } \prod_{\ga \in N(w^{-1})}h_\ga^{mq(w,\ga)}.\ee
\end{theorem}
\noi If we assume hypothesis $(\ref{ii})$ instead of $(\ref{iii})$, it seems difficult to obtain the same estimates on \v Sapovalov elements as in Theorem \ref{Shap}.
However it is still possible to obtain a reasonable estimate using a different definition of the  degree of a partition, at least if $m=1$.
To simplify notation we set $\mathcal{H}_{\gamma} = \mathcal{H}_{\gamma,1}$ and
denote a  \v Sapovalov element for the pair $(\gc, 1)$ by $\gth_{\gamma}$. In the Theorem below we assume that $\Pi$ contains an odd non-isotropic root, since otherwise (\ref{iii}) holds and  the situation is covered by Theorem \ref{Shap}.  This assumption is essential for Lemma \ref{gag}. Likewise
$\gc$ is odd and non-isotropic, then again (\ref{iii}) holds, so we assume that $\gc=w\gb$  with $w \in W_{\rm nonisotropic}$ and
$\gb \in \overline{\Delta}^+_{0} \cup {\overline{\Delta}}^+_{1}$, (see (\ref{Kacroot}) for notation).
\\ \\
 For $\ga$ a positive root, and then for $\pi$ a partition, we define the {\it Clifford degree} of $\ga, \pi$ by $$\cg (\ga) = 2-i, \mbox{ for } \ga \in \Delta^+_{i},\quad   \cg (\pi) = \sum_{\alpha \in \Delta^+} \pi(\alpha)\cg (\ga).$$  The reason for this terminology is that if we set $U_n = \mbox{ span } \{e_{-\pi}|\cg(\pi)\le n\},$ then $\{U_n\}_{n \ge 0}$ is the Clifford filtration on $U(\fn^-)$ as in \cite{M} Section 6.5. The associated graded ring $\gra  U(\fn^-) = \bigoplus_{n \geq 0}U_{n}/U_{n-1}$ is isomorphic to a Clifford algebra. In addition $\gra  U(\fn^-)$ is isomorphic to  an enveloping algebra $U(\fk)$ where the Lie superalgebra $\fk$ is equal to $\fn^-$ as a graded vector space, and the product is modified so that $\fk_0$ is central in $\fk.$
\begin{theorem} \label{aShap} Suppose that $\fg$ is a finite dimensional contragredient, and that $\Pi$ contains an odd non-isotropic root. Assume $\gc$ is a positive root such that $(\ref{i})$ and $(\ref{ii})$ hold.  If $\gc$ is isotropic assume that $m=1.$ Then
there exists a \v Sapovalov  element
$\theta_{\gamma,m} \in U({\mathfrak b}^{-} )^{- m\gamma}$, which is unique modulo the left ideal
$U({\mathfrak b}^{-} ){\mathcal I}(\mathcal H_{\gc})$.  If $m=1$,
the coefficients of $\theta_{\gc}$ satisfy
\be \label{x2}
 2\deg  H_{\pi} \le 2\ell(w)+1 - \cg (\pi),\ee
 and \[H_{\pi^{\gamma}} \mbox{
has leading term } \prod_{\ga \in N(w^{-1})}h_\ga.\]
\end{theorem}
\bc In Theorems \ref{Shap} and  \ref{aShap} $H_{m\pi^{\gamma}}$ is the unique term of highest degree in $\gth_{\gc, m}.$
\ec
\bpf This follows easily from the given degree estimates.\epf
\noi In the next Section we discuss the uniqueness of \v Sapovalov elements. Theorems \ref{Shap} and \ref{aShap} are proved in Sections \ref{s.5} and \ref{s.51} respectively. The proofs depend on a rather subtle cancelation property which is illustrated in Section \ref{cosp}. In Section \ref{s.8} we give a closed formula for \v Sapovalov elements in Type A. 
By definition \v Sapovalov elements give rise  to highest weight vectors in a Verma modules.
The question of when the images of these highest weight vectors in
various factor modules is  non-zero is studied in Section \ref{surv}.
\\ \\
Some of our results hold without assumption (\ref{i}) above.  However it seems more interesting to compare \v Sapovalov elements for an arbitrary Borel to those obtained using the distinguished or anti-distinguished Borel subalgebra as we do in Section \ref{sscbs}. This pair of reference Borels lie at two extremes.  To explain what this means suppose $\fb, \fb'$ is an arbitrary pair of Borel subalgebras (always with the same even part).  Then there is a  sequence
 \be \label{distm} \mathfrak{b} = \mathfrak{b}^{(0)}, \mathfrak{b}^{(1)}, \ldots,
 \mathfrak{b}^{(r)}. \ee
  of Borel subalgebras such that $\mathfrak{b}^{(i-1)}$ and
  $\mathfrak{b}^{(i)}$ are adjacent for $1 \leq i \leq r$, and $\mathfrak{b}^{(r)} = \mathfrak{b}'$.
  If there is no  chain of adjacent Borel subalgebras connecting $\fb$ and $\fb'$ of shorter length than (\ref{distm}), we set $d(\fb,\fb')=r.$ Then $d(\fb,\fb')$ is maximal if and only if $(\fb,\fb')$ is our reference pair, and if $\fb''$ is any other Borel we can
find a chain as in (\ref{distm}) with  $\fb, \fb'$ distinguished and antidistinguished respectively, and $\fb = \mathfrak{b}^{(i)}$ for some $i$ . In addition if $\fg$ is of Type A, B or D and $\gc$ a root of the form $\gd_i-\gep_j$, then $\gc$ (resp. $-\gc$) is a positive root for $\fb$ (resp. $\fb'$) and we can arrange that $\gc$ (resp. $-\gc$) is a simple root for $\mathfrak{b}^{(i-1)}$ (resp. $\mathfrak{b}^{(i)}$).
\brs \label{rmd}{\rm
(a) It is interesting to compare the inequalities (\ref{x1}) and (\ref{x2}).  If $\Pi$ does not contain an odd non-isotropic root, then $\cg(\pi) = 2|\pi|$ for all ${\pi \in {\overline{\bf P}}(\gamma)}$.  Also it  follows by  induction and (\ref{Nw}), that $\hgt \gc = 1 +
\sum_{\ga \in N(w^{-1})} q(w,\ga)$. It follows from (\ref{x1}) that
$$\cg(\pi) + 2\deg  H_{\pi} \le 2+ 2\sum_{\ga \in N(w^{-1})}{q(w,\ga)}.$$
Thus if $q(w,\ga) = 1$ for all  $\ga \in N(w^{-1})$ we have
$$\cg(\pi) + 2\deg  H_{\pi}  \le 2+2\ell(w) .$$
Under the  hypotheses of Theorem \ref{aShap}, (\ref{x2}) sharpens this bound. 
On the other hand if $\fg =\so(5)$ with simple roots $\ga=\gep_1-\gep_2, \gb=\gep_2$, and $\gc = s_\ga(\gb)$, then the inequality (\ref{x2}) does not hold.\\ \\
(b) In general we have
\be \label{gag1} W_{\rm even} \subseteq W_{\rm nonisotropic}\subseteq W.
\ee
Suppose we  use the distinguished set of simple roots.  Then if $\fg$ has Type A, equality holds in \eqref{gag1}, and if $\fg=\osp(1,2n)$ then $W_{\rm even}$ is properly contained in $W_{\rm nonisotropic} = W.$  Otherwise, all
nonisotropic roots contained in the
distinguished  set of simple roots are even thus
$W_{\rm even} = W_{\rm nonisotropic}$. This is however a proper subgroup of $W$.
}\ers
\section{Uniqueness of \v Sapovalov elements.} \label{s.1}
\noi The uniqueness of \v{S}apovalov elements is easily taken care of, and so we do so here.
To do this we need a version of the Jantzen sum formula, which will also play an important role in the sequel.\\ \\
First set \be \label{Kacroot}  {\overline{\Delta}}^+_{0} = \{ \alpha \in \Delta^+_{0} | \alpha
/ 2 \not\in \Delta^+_{1}\}, \quad {\overline{\Delta}}^+_{1}= \{ \alpha
\in \Delta^+_{1}|2\alpha \not\in \Delta^+_{0}\}.\ee
Then for $\lambda \in
\mathfrak{h}^*$ define
\[A(\lambda)_{0}  =  \{ \alpha \in \overline{\Delta}^+_{0} | (\lambda + \rho,
\alpha^\vee) \in \mathbb{N} \backslash \{0\} \}, \]
\[A(\lambda)_{1} = \{ \alpha \in \Delta^+_{1} \backslash \overline{\Delta}^+_{1} |
(\lambda + \rho, \alpha^\vee ) \in 2\mathbb{N} + 1 \}, \]
 \[A(\lambda) = A(\lambda)_{0} \cup  A(\lambda)_{1},\] and
 \[B(\lambda) = \{ \alpha \in \overline{\Delta}^+_{1} | (\lambda + \rho,
\alpha) = 0 \} . \]
If $\alpha \in {\overline{\Delta}}^+_{1},$ let ${\bf p}_{\alpha}(\eta)$ be the number of partitions $\pi$ of $\eta$ such that $\pi{(\ga)} = 0,$ and
then let
$p_{\ga}$ be the generating function given by
$p_\ga =  \sum {\bf p}_\ga(\eta)\epsilon^{-\eta}$.\\ \\
Then by \cite{M} Theorem 10.3.1,  the Jantzen filtration $\{{M}_{i}(\lambda)\}_{i\ge1}$ on ${M}(\lambda)$ satisfies the sum formula
\be \label{jfn}
 \sum_{i > 0} \ch {M}_{i}(\lambda) = \sum_{\alpha \in A(\lambda)}
\ch {M}(s_\alpha \cdot  \lambda) + \sum_{\alpha \in
B(\lambda)} \epsilon^{\lambda -\alpha}p_\ga.\ee
\bt \label{7612} Suppose $\gth_1, \gth_2$ are \v Sapovalov elements for the pair $(\gc,m)$.  Then
\bi \itema for all $\gl \in \mathcal H_{\gc, m}$ we have $\gth_1 v_\gl =\gth_2 v_\gl$
\itemb $\gth_1- \gth_2
\in U({\mathfrak g}){\mathfrak n}^+  + U({\mathfrak g}){\mathcal I}(\mathcal H).
$ \ei\et \bpf
Set
\[\gL = \{\gl \in \mathcal H_{\gc, m}| A(\gl) = \{\gc\}, \;B(\gl) = \emptyset \},\]
if $\gc$ is non-isotropic, and
\[\gL = \{\gl \in \mathcal H_{\gc}| B(\gl) = \{\gc\}, \;A(\gl) = \emptyset \},\]
if $\gc$ is isotropic.  If $\gl \in \gL$
it follows from the sum formula that
$ M_1(\gl)^{\gl-m\gc}$ is one-dimensional.
Because ${M}_{1}(\lambda)$ is the unique maximal submodule of
${M}(\lambda)$, $\gth_1 v_\gl$ and $\gth_2 v_\gl$ are proportional. Then from the requirement that $e_{-\pi^0}$ occurs with coefficient 1 in
a \v{S}apovalov element we have $\gth_1 v_\gl =\gth_2 v_\gl$.  Since $\gL$ is Zariski dense in $\mathcal H_{\gc, m}$, (a) holds and (b) follows from (a) because
by \cite{M} Lemma 9.4.1 we have
\be \label{rut} \bigcap_{\gl \in \Lambda} \ann_{U(\sfg)}v_\gl =
 U({\mathfrak g}){\mathfrak n}^+  + U({\mathfrak g}){\mathcal I}(\mathcal H).
 \ee
\epf
\noi We remark that this proof does not resolve the issue of whether (\ref{kit}) holds in general for Lie superalgebras, but we note that the analog of (\ref{kit}) fails for parabolic Verma modules over simple Lie algebras, \cite{IS}, \cite{IS1}.
 \section{Outline of the Proof and Preliminary Lemmas.}\label{s.4}
Theorems \ref{Shap} and \ref{aShap}
are proved by looking at the proofs given in \cite{H2} or \cite{M} and keeping track of the coefficients. 
Given $\gl \in {\mathfrak h}^* $ we define the {\it specialization} at $\gl$ to be the
map
$$\varepsilon^\gl:U({\mathfrak b}^{-}) = U({\mathfrak n}^-)\otimes S({\mathfrak h})  \longrightarrow M (\lambda),\;  \quad
\sum_{i} a_i \otimes b_i \longrightarrow \sum_{i} a_i  b_i(\gl)v
_{\lambda}.$$
Let  $(\gc,m)$ be as in the statement of the Theorems and set ${\mathcal H} = {\mathcal H}_{\gc, m}$.
If $\gth$ is as in the conclusion of the Theorem, then for any $\gl \in  {\mathcal H},$ $\theta(\gl)v_\gl$
 is a highest weight vector
in $M (\lambda)^{\gl-m\gamma}$. Conversely suppose that $\Lambda$
is a dense subset of ${\mathcal H} $  and that for all $\gl \in
\Lambda$ we have constructed $\theta^\gl \in U({\mathfrak n})^{-m\gamma}$
such that $\theta^\gl v_\gl$ is  a highest weight vector in
$M(\lambda)^{\gl-m\gamma}$ and that
$$\theta^\gl = \sum_{\pi \in {\overline{\bf P}}(m\gc)}a_{\pi, \gl}e_{-\pi}.$$
where $a_{\pi, \gl} $ is a polynomial function of $\gl \in \Lambda$ satisfying suitable conditions.
For $\pi \in {\overline{\bf P}}(m\gc)$,
the assignment $\gl
 \rightarrow a_{\pi, \gl}$ for $\gl \in \Lambda$ determines a polynomial map from ${\mathcal H}$ to
 $U({\mathfrak n}^-)^{-\gamma},$ so  there exists  an element  $H_\pi \in U({\mathfrak h})$ uniquely determined modulo ${\mathcal I}(\mathcal H)$
such that $H_{\pi}(\gl) = a_{\pi, \gl}$ for all $\gl \in \Lambda$.   We define the element
$\theta \in U(\fb^-)$  by
setting
$$\theta= \sum_{\pi \in {\overline{\bf P}}(m\gc)}e_{- \pi}H_{\pi}.$$
\noi Note that $\gth$ is uniquely determined modulo the left ideal
$U({\mathfrak b^-}){\mathcal I}(\mathcal H),$ and that
$\theta(\gl) = \theta^\gl$.
Also, for $\alpha \in \Delta^+$ and
 $\gl \in \Lambda$ we have $e_\alpha \theta v_\gl = e_\alpha
 \theta^\gl v_\gl = 0$, because $\theta^\gl v_\gl = 0$ is a highest weight vector, so  $e_\alpha \theta \in \bigcap_{\gl \in \Lambda} \ann_{U(\sfg)}v_\gl.$   Thus (\ref{boo}) follows from (\ref{rut}).\\ \\
\noi We need to examine the polynomial nature of the coefficients of $\gth_{\gc,m}$. The following easy observation (see \cite{D} Lemma 7.6.9),  is the key to doing this.
Let $A$ be a $\mathbb{Z}_2$-graded associative algebra,
     and suppose that $e$ is an even element of $A$.
     Then for  $a \in A$
     and all $r \in \mathbb{N} $,
\be \label{cow} e^{r}a = \sum^r_{i=0} \left( \begin{array}{c}
                r \\
                i \end{array}\right) ((\ad e)^i a)e^{r - i} .\ee
\noi The following consequence is well-known, \cite{BR}.  We give the short proof for completeness.
\bc \label{oreset} With the same hypothesis as above, suppose that $\ad e$ is locally nilpotent.   Then the set $\{e^n|n \in \mathbb{N}\} $  is an Ore set in $A$. \ec \bpf Given $a \in A$  and $n \in \mathbb{N}$, suppose that
$(\ad e)^{k+1} a = 0$.  Then $e^{k+n}a = a'e^n$, where \[a' = \sum^k_{i=0} \left( \begin{array}{c}
                n+k \\
                i \end{array}\right)
                ((\ad e)^i a)e^{k-i}.\]
 \epf
\noi Now suppose $\ga \in \Pi_{\rm nonisotropic}$, and set $e = e_{-\ga}.$ Then $e$ is a nonzero divisor in
$U = U(\fn^-)$, and the set $\{e^n|n \in \mathbb{N}\} $ is an Ore set in $U $ by Corollary \ref{oreset}. We write $U_e$ for the corresponding Ore localization. The adjoint action of $\fh$ on $U$ extends to $U_e$, and in the next result we give a basis for the weight spaces of $U_e$.
Let ${\widehat{\bf P}}(\eta)$ be the set of pairs $(k,\pi)$ such that $k \in \mathbb{Z}, \pi \in {\overline{\bf P}}(\eta -k\ga)$ and $\pi(\ga) = 0.$  Then we have \bl \label{uebasis}  \begin{itemize} \item[{}]
\itema The set $\{e_{-\pi} e^k|(k,\pi ) \in {\widehat{\bf P}}(\eta)\}$ forms a $\ttk $-basis for the weight space $U_e^{-\eta}.$
 \itemb If $u = \sum_{(k,\pi ) \in {\widehat{\bf P}}(\eta)} c_{(k,\pi )} e_{-\pi} e^k \in U_e^{-\eta}$ with $c_{(k,\pi )} \in \ttk ,$ then $u \in U$ if and only if $c_{(k,\pi )} \neq 0$ implies $k \ge 0$.
\end{itemize}
 \el \bpf (a) Suppose $u \in U_e^{-\eta}$. We need to show that $u$ is uniquely expressible in the form
 \be \label{china} u = \sum_{(k,\pi ) \in {\widehat{\bf P}}(\eta)} c_{(k,\pi )} e_{-\pi} e^k \ee
 We have $ue^N \in U^{-(N\ga+\eta)}$ for some $N$.  Hence by the PBW Theorem for $U$ we have a unique expression $$ue^N = \sum_{\gs \in {\bf\overline{P}}(\eta+N\ga)}a_\gs  e_{-\gs}  .$$ Now if $a_\gs \neq 0$, then $e_{-\gs}= e_{-\pi}e^{\ell}$ where $\gs(\ga) = {\ell}$ and
  where $\pi \in {\bf\overline{P}} (\eta+N\ga-{\ell}\ga)$ satisfies $\pi(\ga) = 0$. Then $\pi$ and ${k = N-\ell}$ are uniquely determined by $\gs$, so we set $c_{(k,\pi )} =a_\gs$.  Then clearly (\ref{china})  holds. Given (a), (b) follows from the PBW Theorem.
  \epf
\noi  We remark that if $\ga$ is a non-isotropic odd root, then we can use $e^2$ in place of $e = e_{-\ga}$ in the above Corollary and Lemma.  \noi   However we will need a version of Equation (\ref{cow}) when $e$ is replaced by an odd element $x$ of a $\Z_2$-graded algebra $A$.
Suppose that $z$ is homogeneous, and define $z^{[j]} = (\ad x)^j z$.  Set $e = x^2$ and apply (\ref{cow}) to $a = xz =[x, z] + (-1)^{\overline z}zx,$  to obtain
\be \label{f11}  x^{2\ell+1}z = \sum^\ell _{i=0} \left( \begin{array}{c}
                \ell \\
                i \end{array}\right) z^{[2i+1]}x^{2(\ell - i)}
 + (-1)^{\overline z}\sum^\ell _{i=0} \left( \begin{array}{c}
                \ell \\
i \end{array}\right) z^{[2i]} x^{2\ell - 2i +1} .\ee

\noi The \v Sapovalov elements in Theorems \ref{Shap} and  \ref{aShap} are constructed inductively using the next Lemma.
Suppose that the pair $(\gc, m)$ satisfies one of the following
\be \label{fat} m \mbox{ is an odd positive integer if } {\gc} \mbox{ is an odd non-isotropic root,}\ee
\be \label{bog} m \mbox{ is  a positive integer if } {\gc} \mbox{ is an even root such that } {\gc}/2 \mbox{ is not a root.} \ee
\be \label{log}{\gc} \mbox{ is an odd isotropic root and } m = 1.\ee
\begin{lemma}\label{1768}
Suppose that  $\mu \in \cH_{\gc',m}, \alpha \in
\Pi_{\rm nonisotropic}, $ and set 
\be \label{21c}\mu = s_\alpha\cdot \gl,\;\gc' = s_\alpha\gc,\;
  p = (\mu + \rho, \alpha^\vee),\;
q = (\gc, \alpha^\vee).\ee
Assume that $p\ge 0, q  \in \mathbb{N}
\backslash \{0\}$,
the pair $(\gc,m)$ satisfies {\rm(\ref{fat})} or {\rm(\ref{bog})},
and
\begin{itemize}
\itema  $\theta' \in
U({\mathfrak n}^-)^{-m\gc'}$ is such that $v = \theta'v_\mu \in
M(\mu)$ is a highest weight vector.
\itemb If $\alpha \in \overline{\Delta}^+_{0}$, then $q=1.$
\itemc  If  $\alpha \in \Delta^+_{1} \backslash \overline{\Delta}^+_{1} $ then $q = 2$ and $p$ is odd.
\end{itemize}
\noi Then there is a unique $\theta \in U({\mathfrak n}^-)^{-m\gc}$ such that
\begin{equation} \label{21nd}
e^{p + mq}_{- \alpha}\theta' = \theta e^p_{- \alpha}.
\end{equation}
\end{lemma}
\bpf This is well-known, see for example  \cite{H2} Section 4.13 or \cite{M} Theorem 9.4.3.\epf 
\noi Equation \eqref{21nd} is the basis for the proof of many properties of \v Sapovalov elements. We note the following variations. First under the hypothesis of the Lemma the \v Sapovalov elements
$\theta_{\gamma',m}$ and 
$\theta_{\gamma,m}$ are related by 
\begin{equation} \label{21a}
e^{p + mq}_{- \alpha}\theta_{\gc',m}(\mu) = \theta_{\gc,m}(\gl) e^p_{- \alpha}.
\end{equation}
Now set $r  = (\gl + \rho-m\gc, \alpha^\vee)$ and instead of the hypothesis on $p$, suppose that $r
\ge 0$.  Then 
\begin{equation} \label{21b}
\theta_{\gc',m}(\mu)
e^{r + mq}_{- \alpha} = e^r_{- \alpha} \theta_{\gc,m}(\gl).
\end{equation} Note that \eqref{21b} becomes  formally equivalent to \eqref{21a} when we set 
$r=-{(p + mq)}$.\\ \\
\noi In the proofs of Theorems \ref{Shap} and \ref{aShap} we write $\gc = w\gb$ for $\gb \in \Pi$ and $w\in  W$.
We use the Zariski dense  subset $\gL$ of  ${\mathcal H}_{\beta,m}$ defined by
\be \label{tar}\Lambda =
\left\{
\begin{array}{c|c}
\nu \in {\mathfrak h}^* & (\nu + \rho, \beta^\vee)=m \mbox{ and }(\nu + \rho,
\alpha^\vee) \in \mathbb{N} \backslash
\{0 \}\mbox{ for}\\ &  \mbox{all }  \alpha \in \Pi, \mbox{ with }(\nu + \rho,
\alpha^\vee)
\mbox{ odd if }
\alpha
\mbox{
is odd}
\end{array}
\right\}
.\ee
\section{Proof of Theorem \ref{Shap}.}\label{s.5}
In this section we assume $\fg$ is contragredient and hypotheses $(\ref{i})$ and $(\ref{iii})$ hold.  If $\gc$ is a simple root, then $\theta_{\gamma,m} = e_{-\gc}^m$ satisfies the conditions of Theorem \ref{Shap}.
Otherwise we have $\gc = w\gb$ for some  $w \in W_{\rm even}, w
\neq 1$. 
Write \be \label{wuga} w =
s_{\ga}u, \quad   \gamma' = u\gb, \quad
\gamma = w\gb =
s_{\alpha} \gamma', \ee
 with $\alpha \in \Pi_{\rm even}$ and $\ell(w) =
\ell(u) + 1.$ Since the statement of Theorem \ref{Shap} involves precise but somewhat lengthy
conditions on the coefficients, we introduce the following definition as a shorthand.
\bd \label{sag} We say that a family of elements $\theta^\gl_{\gamma,m} \in
U({\mathfrak n}^-)^{- m\gamma}$ is {\it well posed for $w$} if  for all
$\gl \in w\cdot  \Lambda$ we have \ed
\begin{equation} \label{2}
\theta^\gl_{\gamma,m} = \sum_{\pi \in {\overline{\bf P}}(m\gamma)} a_{\pi,\gl}
e_{-\pi},
\end{equation}
where
the coefficients $a_{\pi,\gl} \in \ttk $ depend polynomially on $\gl \in
w\cdot \Lambda$, and \bi%
\itema $\deg  \; a_{\pi,\gl} \leq m\hgt \gc - |\pi|$
\itemb $a_{m\pi^{\gamma},\gl}$ is a polynomial function of
$\gl$ of degree $m(\hgt \gc - 1)$ with highest term equal to $c\prod_{\ga
\in N(w^{-1})}(\gl,\ga)^{mq(w,\ga)}$ for a nonzero constant $c$.
\end{itemize}
We show that the conditions on the coefficients in this definition
are independent of the ordering on the positive roots $\Delta^+$ used
to define the $e_{-\pi}.$
  Consider  two orderings on $\Delta^+$, and for $\pi \in
\overline{\bf P}(m\gc)$, set $e_{-\pi} = \prod_{\alpha \in \Delta^+} e^{\pi
(\alpha)}_{- \alpha}$, and $\overline{e}_{-\pi} = \prod_{\alpha \in
\Delta^+} e^{\pi (\alpha)}_{- \alpha}$, the product being taken with
respect to the given orderings. \bl \label{delorder1}%
Fix a total order on the set ${\overline{\bf P}}(m\gc)$ such that if $\pi, \gs \in {\overline{\bf P}}(m\gc)$ and $|\pi| > |\gs|$ then $\pi$ precedes $\gs,$ and use this order on partitions to induce orders on the bases ${\bf B_1} = \{
e_{-\pi}|\pi \in {\overline{\bf P}}(m\gc)\}$  and ${\bf B_2} = \{\overline{e}_{-\pi}|\pi \in {\overline{\bf P}}(m\gc)\}$  for $U(\fn^-)^{-m\gc}$. Then the change of basis matrix from the basis ${\bf B_1}$ to ${\bf B_2}$ is upper triangular with all diagonal entries equal to $1.$
\el %
\noindent \bpf Let $\{ U_n =U_n(\fn^-) \}$ be the standard filtration on
$U=U({\mathfrak n}^-).$  Note that if $\pi \in {\overline{\bf P}}(m\gc),$ then $e_{-\pi}, \overline{e}_{-\pi} \in
U_{|\pi|}({\mathfrak n}^-)^{-m\gc}.$  Also the factors of $e_{-\pi}$ commute
modulo lower degree terms, so for all $\pi \in {\overline{\bf P}}(m\gc),$ $e_{-\pi} -\overline{e}_{-\pi} \in
U_{|\pi| - 1}({\mathfrak n}^-)^{-m\gc}$.  The result follows easily.
 \epf
\bl \label{delorder}%
For $x \in U({\mathfrak n}^-)^{-m\gc}
\otimes S({\mathfrak h}),$ write
\be \label{owl} x \; = \; \sum_{\pi \in {\overline{\bf P}}(m\gc)}
e_{-\pi}f_\pi \; = \; \sum_{\pi \in {\overline{\bf P}}(m\gc)}
\overline{e}_{-\pi}g_\pi.\ee Suppose that $f_{m\pi^{\gamma}}$ has
degree $m(\hgt \gc - 1),$ and that for all $\pi \in {\overline{\bf P}}(m\gc),$ we have $\deg  f_\pi \le m\hgt \gc - |\pi|.$  Then
$g_{m\pi^{\gamma}}$ has the same degree and leading term as
$f_{m\pi^{\gamma}}$
 and
 for all $\pi \in {\overline{\bf P}}(m\gc),$ we
have $\deg  \;g_\pi \le m\hgt \gc - |\pi|.$
\el %
\noindent \bpf By Lemma \ref{delorder1} we can write
\[e_{-\pi} = \sum_{\zeta
 \in {\overline{\bf P}}(m\gc)}
c_{\pi,\zeta}\overline{e}_{-\zeta},\] where $c_{\pi,\zeta} \in \ttk ,
c_{\pi,\pi} = 1$ and if $c_{\pi,\zeta} \neq 0$ with $\zeta \neq
\pi,$ then $|\gz| < |\pi|$. Thus (\ref{owl}) holds with
$$g_\zeta = \sum_{\pi \in {\overline{\bf P}}(m\gc)}
c_{\pi,\zeta}f_{\pi}.$$ It follows that $g_\zeta$ is a linear combination of polynomials of
degree less than $m\hgt \gc - |\gz|.$ Also $|m\pi^{\gamma}| = m,$ and for $\zeta \in {\overline{\bf P}}(m\gc), \zeta \neq
m\pi^{\gamma},$ we have $|\gz|>m.$  Therefore $$g_{m\pi^{\gamma}}
= f_{m\pi^{\gamma}} + \mbox{ a linear combination of polynomials of smaller
degree.}$$ The result follows easily from this.\epf
\noi Now recall the notation from Equation (\ref{wuga}).  Suppose $\nu \in \Lambda$ and set  \be \label{wuga1}
\mu =
u\cdot \nu, \quad \lambda = w\cdot \nu
= s_{\alpha}\cdot \mu.
 \ee
\noi The next Lemma is the key step in establishing the degree estimates in the proof of Theorem \ref{Shap}.
The idea is to use Equation (\ref{21nd}) and the fact that
$\theta \in U({\mathfrak n}^-)$, rather than a localization of $U({\mathfrak n}^-)$, to show that certain coefficients cancel.
Then using induction and (\ref{21nd}) we obtain the required degree estimates. Since the proof of the Lemma is rather long we break it into
a number of steps.
\bl \label{wpfg} Suppose that $p, m, q$ are as in Lemma \ref{1768}, $\ga \in \Pi_{\rm even}$ and
\begin{equation} \label{2nd} e^{p+mq}_{-
\alpha}\theta^{\mu}_{\gamma',m} = \theta^\lambda_{\gamma,m}e^p_{-
\alpha}.
\end{equation} 
Then the family $\theta^\gl_{\gamma,m} $ is well posed for $w$ if
the family $\theta^\mu_{\gamma',m}$ is well posed for $u$. \el \bpf
{\it Step 1. Setting the stage.}\\  \\
 Suppose that \begin{equation} \label{4cha7} \theta^{\mu}_{\gamma',m} =
\sum_{\pi' \in {\overline{\bf P}}(m\gamma')} a'_{\pi',\mu} e_{-\pi'},
\end{equation}
and let
\be \label{qwer} e_{- \pi'}^{(j)} = (\ad e_{- \alpha})^j e_{- \pi'} \in U_{|\pi'|}(\fn^-)^{-(m\gc' +
j \ga )},\ee
for all $j \geq 0,$ and $\pi' \in {\overline{\bf P}}(m\gamma')$. Then by Equation (\ref{cow})
\begin{eqnarray} \label{49}
e^{p + mq}_{- \alpha} e_{- \pi'} = \sum_{i \geq 0} \left(
\begin{array}{c}
p + mq \\
j
\end{array}
\right) e^{(j)}_{- \pi'} \; e^{p + mq-j}_{- \alpha}.
\end{eqnarray}
Choose $N$ so that
$e^{(N+1)}_{- \pi'} = 0,$ for all $\pi' \in {\overline{\bf P}}(m\gc').$ Then for
all such $\pi' $ and $j = 0, \ldots, N$ we can write
\begin{eqnarray} \label{99}
e_{- \pi'}^{(j)}e^{N-j}_{- \alpha} = \sum_{\zeta \in {\overline{\bf P}}(m\gc' +
N \ga )} b_{j, \zeta}^{\pi'}e_{- \zeta} ,
\end{eqnarray}
with $ b_{j, \zeta}^{\pi'} \in \ttk .$ Furthermore if  $ b_{j,
\zeta}^{\pi'} \neq 0,$ then
since
$e_{- \pi'}^{(j)}e^{N-j}_{- \alpha} \in  U_{|\pi'| + N - j}$,
 (\ref{qwer}) gives
\begin{eqnarray} \label{98}
|\gz| \leq |\pi'| + N - j.
\end{eqnarray}
{\it Step 2. The cancelation step.}\\ \\
By Equations (\ref{4cha7}) and (\ref{49})
\begin{eqnarray} \label{199} e^{p + mq}_{- \alpha}
\theta^{\mu}_{\gamma',m}  & = & \sum _{\pi' \in
{\overline{\bf P}}(m\gamma')} a'_{\pi',\mu} 
e^{p + mq}_{- \alpha}e_{-\pi'}\\
& = & \sum_{ \pi' \in {\overline{\bf P}}(m\gamma'),\;j \geq 0} \left(
\begin{array}{c}
p + mq \nonumber\\
j
\end{array}
\right) a'_{\pi',\mu} e_{- \pi'}^{(j)}e^{p + mq-j}_{- \alpha}.
\end{eqnarray}
Now collecting coefficients, set \be\label{59} c_{ \zeta,\lambda}
 = \sum_{\pi' \in {\overline{\bf P}}(m\gamma'),\;j \geq 0} \left(
\begin{array}{c}
p + mq \\
j
\end{array}
\right) a'_{\pi',\mu} 
b_{j, \zeta}^{\pi'}. \ee Then using Equations (\ref{99}) and (\ref{199}), we have in $U_e,$ where $e = e_{-\ga}$,
\begin{eqnarray} \label{lab} e^{p + mq}_{- \alpha} \theta^{\mu}_{\gamma',m} &
= & \sum_{\zeta \in {\overline{\bf P}}(m\gc' + N\ga )}%
c_{ \zeta, \lambda}e_{-\zeta}e^{p + mq-N}_{- \alpha}.
\end{eqnarray}
By (\ref{2nd}) and Lemma \ref{uebasis}, $c_{ \zeta, \lambda} = 0$ unless $\zeta(\ga) \geq
N-mq.$\\ \\
{\it Step 3. The coefficients $a_{\pi,\lambda}$.}\\ \\
It remains to deal with the nonzero terms $c_{ \zeta, \lambda}$. There is a bijection
\[f:{\overline{\bf P}}(m\gamma) \longrightarrow \{\zeta \in {\overline{\bf P}}(m\gamma' +N\ga)\;|\;\zeta(\ga) \geq N-mq\},\]
defined by \be \label{yam}
 (f\pi)(\gs)=\left\{ \begin{array}
  {ccl}\pi(\gs)&\mbox{if} \;\; \gs \neq \ga,
\\\pi(\ga) + N - mq &\mbox{if} \;\; \gs = \ga.
\end{array} \right. \ee
Moreover if $f\pi = \zeta,$ then
\begin{eqnarray} \label{119}|\gz| = |\pi| + N -mq
\end{eqnarray} and $e_{-\pi} = e_{-\zeta}e^{mq-N}_{-
\alpha}$.
Thus in Equation (\ref{2nd}) the coefficients $a_{\pi,\lambda}$
of $\theta^\lambda_{\gamma,m} $ (see (\ref{2})) are given by \be
\label{cwp2} a_{\pi,\lambda} = c_{f(\pi),\lambda}.\ee
{\it Step 4. Completion of the proof.}\\ \\
We now show that the family $\theta^\gl_{\gamma,m} $ is well posed for $w.$ For this we
use Equations (\ref{59}) and  (\ref{cwp2}), noting that $p =
(s_\alpha\cdot \gl + \rho,\ga)$ depends linearly on $\gl.$   It is clear
that the coefficients $a_{\pi,\lambda}$ are polynomials in
$\lambda$.
 By induction $\deg  a'_{\pi',\mu} \le m\hgt \gc' - |\pi'|.$ Thus using
 (\ref{59}), \be \label{fox} \; \mbox{deg} \;\; a_{\pi,\gl} = \;
\mbox{deg} \; c_{\zeta,\gl}  \leq  \max\{j + \; \mbox{deg} \;
a'_{\pi',\mu}\;|\; b_{j, \zeta}^{\pi'}
\neq 0\}. \ee Now if $b_{j, \zeta}^{\pi'}
\neq 0$ then Equation (\ref{98}) holds. Therefore by Equation (\ref{119}) 
\begin{eqnarray} 
\mbox{deg} \;\; a_{\pi,\gl} \nonumber &\leq&
\mbox{deg} \;\; a'_{\pi',\mu}+ |\pi'| + N - |\gz|
\\ &=&  \mbox{deg} \;\; a'_{\pi',\mu}+ |\pi'| - |\pi| +mq\nonumber 
\end{eqnarray}
Finally since $\gamma = \gc' +
q\ga$, induction gives
(a) in Definition \ref{sag}.\\ \\
\noi  Also, modulo terms of lower
degree
\be \label{Nw1}  a_{m\pi^{\gamma},\gl} \equiv \left( \begin{array}{c}
p + mq\\
mq
\end{array}   \right)
a'_{m\pi_{\gamma'},\mu} .\ee
Note that the above binomial coefficient is a polynomial of degree $mq$ in $p$.
 By induction 
 $a'_{m\pi_{\gamma'},\mu}$ has
highest term $c'\prod_{\gt \in N(u^{-1})}(\mu,\gt)^{mq(u,\tau)}$ as a polynomial in $\mu$, for
a nonzero constant $c'$. Now $(\mu + \gr,\gt) = (\gl +\gr, s_\ga \gt)$ and $(\mu + \gr,\gt)- (\mu,\gt)$, $(\gl +\gr, s_\ga \gt) - (\gl, s_\ga \gt)$ are
constant.  Therefore  as a polynomial in $a'_{m\pi_{\gamma'},\mu}$ has
highest term $$c'\prod_{\gt \in N(u^{-1})}(\gl,s_\ga\gt)^{mq(u,\tau)} = c'\prod_{\gt \in N(u^{-1})}(\gl,s_\ga\gt)^{mq(w,s_\ga\tau)}.$$
From the representation theory of $\fsl(2)$, it follows that $(\ad e_{- \alpha})^q(\fg^{-\gc'}) = \fg^{-\gc}$.  Since $e_{-\gc}$ is not used in the construction of $\theta_{\gamma',m},$ we can choose the notation so that
$e^{(mq)}_{- m\gamma'} = e_{- m\gamma}$. Then $e^{(mq+1)}_{- m\gamma'} = 0.$ Now $q(w,\ga) = (\gc, \alpha^\vee) = q,$ and the degree of the binomial coefficient in (\ref{Nw1}) as a polynomial in $p$ is $mq$, so the
claim about the leading term of $a_{m\pi^{\gamma},\gl} $ in Definition \ref{sag} (b) follows from
Equations (\ref{Nw1}) and (\ref{Nw}).\epf
\begin{theorem}  \label{localshap}
Suppose $\gamma = w
\beta$ with $w \in W_{\rm even}$ and $\gb$ simple. There exists a family of elements
$\theta^\gl_{\gamma,m} \in U({\mathfrak n}^-)^{- m\gamma}$ for all $\gl \in
w\cdot  \Lambda$  which is well posed for $w$ such that
\[\theta^\gl_{\gamma,m}v_{\gl} \; \mbox{is a highest weight vector in}\;
M(\gl)^{\gl - m\gamma}.\]
\end{theorem}
\noindent \bpf We use induction on the length of $w.$ If $w = 1,$ we
take $\theta^\gl_{\gamma,m} = e_{-\gb}^m$ for all $\gl.$
Now assume that $w
\neq 1,$ and use the notation of Equations (\ref{wuga}) and (\ref{wuga1}).
If $\lambda =
 s_{\alpha}\cdot \mu \in \Lambda,$ then it is well known that
$M(\gl)$ is uniquely embedded in
 $M(\mu).$
 Set \[p = (\mu + \rho,
\alpha^{\vee}) = (\nu + \rho, u^{-1}\alpha^{\vee}), \quad
 (\gamma, \alpha^{\vee}) = q.\] Then $p$ and $q$ are positive
integers. Also $\gl = \mu - p\alpha$ and $\gamma = \gamma' +
q\alpha.$
By induction
there exist  elements $\theta^{\mu}_{\gamma',m} \in U({\mathfrak n}^-)^{-
\gamma'}$ which are well posed for $u$ such that
\[v = \theta^{\mu}_{\gamma',m}v_{\mu} \in M(\mu)^{\mu -
m\gamma'} \mbox{is a highest weight vector}.
\]
By Lemma \ref{1768} there exists a unique element
$\theta^{\lambda}_{\gamma,m} \in U({\mathfrak n}^-)^{- m\gamma}$ such that
(\ref{2nd})  holds
and therefore
\[e^{p+mq}_{- \alpha} v
=
\theta^{\lambda}_{\gamma,m}
e^{p}_{- \alpha} v_{\mu}
\in U({\mathfrak n}^-)e^p_{- \alpha} v_{\mu} =
M(\lambda).\]
It follows from Lemma \ref{wpfg} that the family
$\theta^\gl_{\gamma,m}$ is well posed for $w$. \epf \noindent {\it Proof of Theorem \ref{Shap}.} Let $\theta^\gl_{\gamma,m}$ be the family of elements from Theorem \ref{localshap}. The existence of the elements
$$\theta_{\gc, m} = \sum_{\pi \in {\overline{\bf P}}(m\gc)} e_{-\pi} H_{\pi} \in U(\fb^-),$$ such that $\theta_{\gc, m}(\gl) = \theta_{\gc, m}^\gl$ for all $ \lambda \in \Lambda$
follows since $ w \cdot\Lambda$ is
Zariski dense in ${\mathcal H}_{\gc,m}.$
The claims about the coefficients
$H_{\pi}$ follow from the fact that the family $\theta^\gl_{\gamma,m}$ is well posed for $w$.
\hfill  $\Box$
\section{Proof of Theorem \ref{aShap}.} \label{s.51}
Now suppose that $\fg$ is contragredient, with $\Pi$ as in $(\ref{i})$. We assume that
$\Pi$ contains an odd non-isotropic root and $\gc=w\gb$  with $w \in W_{\rm nonisotropic}$,
$\gb \in \overline{\Delta}^+_{0}\cup {\overline{\Delta}}^+_{1}$.
Our assumptions have the following consequence.
\bl \label{gag} Suppose that  $\gc$ is an isotropic root and $\ga \in \Pi_{\rm nonisotropic}$ is such that $(\gc,\ga)\neq 0.$ Then
\bi \itema If $\ga$ is even then $(\gc,\ga^\vee) = \pm 1$.
\itemb  If $\ga$ is odd then $(\gc,\ga^\vee) = \pm 2$. \ei\el \bpf Left to the reader. \epf
\noi In this section $\{ U_n =U_n(\fn^-) \}$ is the Clifford filtration on
$U=U({\mathfrak n}^-).$
\bl \label{it} The Clifford filtration on $U(\fn^-)$  stable under the adjoint action of $\fn_0^-,$ and satisfies $\ad \fn_1^- (U_n) \subseteq U_{n+1}$. \el
\bpf Left to the reader.\epf\noi
\noi Fix a total order on the set ${\overline{\bf P}}(\gc)$ such that if $\pi, \gs \in {\overline{\bf P}}(\gc)$ and $|\pi| > |\gs|,$ or if
$|\pi| = |\gs|$ and $\cg(\pi) > \cg(\gs)$, then $\pi$ precedes $\gs,$ and use this order to induce orders on the bases ${\bf B_1} = \{
e_{-\pi}|\pi \in {\overline{\bf P}}(\gc)\}$  and ${\bf B_2} = \{\overline{e}_{-\pi}|\pi \in {\overline{\bf P}}(\gc)\}$  for $U(\fn^-)^{-\gc}$.
  Consider  two orderings on $\Delta^+$, and for $\pi \in
\overline{\bf P}(\gc)$, set $e_{-\pi} = \prod_{\alpha \in \Delta^+} e^{\pi
(\alpha)}_{- \alpha}$, and $\overline{e}_{-\pi} = \prod_{\alpha \in
\Delta^+} e^{\pi (\alpha)}_{- \alpha}$, the product being taken with
respect to the given orderings.
\\ \\
\noi Next we prove a Lemma relating $|\pi|$ to $\cg(\pi)$.
Set
$$\Xi = {\Delta}^+_{0} \backslash \overline{\Delta}^+_{0}, \mbox{ and }a(\pi)  = \sum_{2\gd \in \Xi}\pi(\gd).$$
\bl Suppose $\pi, \gs \in {\overline{\bf P}}(\gc)$.
\bi \label{pig1}\itema We have   $2|\pi| -\cg(\pi) = a(\pi)$.\ei
\bi \itemb $a(\pi) \le 2.$
\itemc If $\gs$ precedes $\pi$, then $\cg(\pi) \le \cg(\gs).$
\ei \el
\bpf (a) follows since 
\[ |\pi| =a(\pi)+ \sum_{\ga \in \overline{\Delta}_{0}^+}\pi(\ga) ,\]
and \[ \cg(\pi) =a(\pi)+ 2\sum_{\ga \in \overline{\Delta}_{0}^+}\pi(\ga).\]
For  (b) we note that the Lie superalgebras that have 
 an odd non-isotropic root $\gd$ are $G(3)$ and the family $\osp(2m+1,2n)$.
     Define a group homomorphism $f:\bigoplus_{\ga \in \Pi}\Z\ga \lra \mathbb{Z}$ by setting  $f( \gd) = 1$ and $f(\ga) = 0$ for any $\ga \in \Pi$, $\ga \neq \gd.$
It can be checked on a case-by-case basis that if $\gd\in \Pi$, then $\gd$ occurs with coefficient at most two when a positive root $\gc$ is written as a linear combination of simple roots. Since $a(\pi)  = f(\gc)$  for $\pi\in  {\overline{\bf P}}(\gc)$,
(b) follows. If (c) is false, then  by definition of the order, we must have $|\pi| < |\gs|$ and
$\cg(\pi) > \cg(\gs).$ But then by (a) this implies that $$a(\gs) =2|\gs| -\cg(\gs)  \ge 2|\pi| -\cg(\pi)+3 = a(\pi) + 3\ge 3$$ which contradicts (b). \epf
\bd \label{sog} We say that a family of elements $\theta^\gl_{\gamma} \in
U({\mathfrak n}^-)^{- \gamma}$ is {\it well posed for $w$} if  for all
$\gl \in w\cdot  \Lambda$ we have
\begin{equation} \label{2os}
\theta^\gl_{\gamma} = \sum_{\pi \in {\overline{\bf P}}(\gamma)} a_{\pi,\gl}
e_{-\pi},
\end{equation}
where
the coefficients $a_{\pi,\gl} \in \ttk $ depend polynomially on $\gl \in
w\cdot \Lambda$, and \bi%
\itema $2\deg  \; a_{\pi,\gl} \leq 2\ell(w)  +1 - \cg(\pi)$
\itemb $a_{\pi^{\gamma},\gl}$ is a polynomial function of
$\gl$ of degree $\ell(w)$ with highest term equal to $c\prod_{\ga
\in N(w^{-1})}(\gl,\ga)^{q(w,\ga)}$ for a nonzero constant $c$.
\end{itemize}\ed
\bl \label{ant}%
Write $x \in U({\mathfrak n}^-)^{-\gc}
\otimes S({\mathfrak h}),$
\[x \; = \; \sum_{\pi \in {\overline{\bf P}}(m\gc)}e_{-\pi}f_\pi \; = \; \sum_{\pi \in {\overline{\bf P}}(m\gc)}\overline{e}_{-\pi}g_\pi.\]
as in Equation {\rm (\ref{owl})}. If the coefficients $f_\pi,$ satisfy
\bi%
\itema $2\deg  \; f_{\pi} \leq 2\ell(w)  +1 - \cg(\pi)$
\itemb $f_{\pi^{\gamma}}$ is a polynomial  of degree $\ell(w)$ with highest term equal to $c\prod_{\ga
\in N(w^{-1})}h_\ga^{q(w,\ga)}$ for a nonzero constant $c$.
\end{itemize}
then the  coefficients  $g_\pi$ satisfy the same conditions.\el \bpf Taking Lemma \ref{pig1} into account, this is the same as the proof of Lemma {\rm \ref{delorder}}.\epf
\bl \label{wpfg1} Suppose that $p, q$ and $\ga$ are as in Lemma \ref{1768} and
\begin{equation} \label{3nd} e^{p+q}_{-
\alpha}\theta^{\mu}_{\gamma'} = \theta^\lambda_{\gamma}e^p_{-
\alpha}.
\end{equation} 
Then the family $\theta^\gl_{\gamma} $ is well posed for $w$ if
the family $\theta^\mu_{\gamma'}$ is well posed for $u$. \el \bpf
If $\ga$ is odd, then  $p = 2\ell -1$ is odd by (\ref{fat}), and $q =2$ by Lemma \ref{gag}.
Write $\theta^{\mu}_{\gamma'}$
as in (\ref{4cha7}) and then define
the $e_{- \pi'}^{(j)}$ as in (\ref{qwer}).
Set $\gve(\gc') = 1$ if $\gc'$ is an even root and $\gve(\gc') = -1$ if $\gc'$ is odd. Then instead of (\ref{49}) we have, by (\ref{f11})
\be \label{f12}  e_{-\ga}^{2\ell+1}e_{- \pi'} = \sum^\ell _{i=0} \left( \begin{array}{c}
                \ell \\
                i \end{array}\right) e_{- \pi'}^{(2i+1)}e_{-\ga}^{2(\ell - i)}
 + \gve(\gc')\sum^\ell _{i=0} \left( \begin{array}{c}
                \ell \\
i \end{array}\right) e_{- \pi'}^{(2i)} e_{-\ga}^{2\ell - 2i +1} .\ee
\noi Parallel to the definition of  the $b_{j, \zeta}^{\pi'}$ in (\ref{99}), we set for sufficiently large $N$
\[e_{- \pi'}^{(j)}e^{N-j}_{- \alpha}= \sum_{\gz \in {\overline{\bf P}}(\gamma'+N\ga)} b_{j, \zeta}^{\pi'}e_{- \gz}.\]
For $x \in \mathbb{R}$ we denote the largest integer not greater than $x$ by $\left\lfloor x \right\rfloor$. Then if  $ b_{j,
\zeta}^{\pi'} \neq 0,$ we have
\begin{eqnarray} \label{98o}
\cg(\zeta) \leq \cg(\pi') + N - 2\left\lfloor \frac{j}{2} \right\rfloor.
\end{eqnarray}
Indeed this holds because by Lemma \ref{it}, we have for such $j$
\be \label{la1} e_{- \pi'}^{(j)}e^{N-j}_{- \alpha} \in  U_{\cg(\pi') + N - 2\left\lfloor \frac{j}{2} \right\rfloor}.\ee
 Replacing (\ref{59}) we set,
\be\label{59o} c_{ \zeta,\lambda}
= \sum_{\pi' \in {\overline{\bf P}}(\gamma'),\;i \geq 0} \left(
\begin{array}{c}
\ell \\i\end{array}
\right) a'_{\pi',\mu} b_{2i+1, \zeta}^{\pi'}
+ \gve(\gc')\sum_{\pi' \in {\overline{\bf P}}(\gamma'),\;i \geq 0} \left(
\begin{array}{c}
\ell \\i\end{array}
\right) a'_{\pi',\mu} b_{2i, \zeta}^{\pi'}. \ee
Then we obtain the following variant of Equation (\ref{lab})
\begin{eqnarray} \label{labos} e^{2\ell+1}_{- \alpha} \theta^{\mu}_{\gamma'} &
= & \sum_{\zeta \in {\overline{\bf P}}(m\gc + N\ga )}
c_{ \zeta, \lambda}e_{-\zeta}e^{2\ell+1-N}_{- \alpha}.
\end{eqnarray}
In the cancelation step we find that $c_{ \zeta,\lambda}=0$ unless $\gz(\ga)\ge N-2$, and the bijection
\[f:{\overline{\bf P}}(\gamma) \longrightarrow \{\zeta \in {\overline{\bf P}}(\gamma' +N\ga)\;|\;\zeta(\ga) \geq N-2\},\]
is defined as in Equation (\ref{yam}) with $m=1$ and $q=2.$ Then the coefficients  $a_{\pi,\lambda}$ are defined as in (\ref{cwp2}).
Instead of Equations (\ref{119})  and  (\ref{fox}) we have
\begin{eqnarray} \label{119o}\cg(\zeta) = \cg(\pi) + N -2,
\end{eqnarray}
and \be \label{foxa} \; \mbox{deg} \;\; a_{\pi,\gl} \;
 \leq  \max\{\lfloor j/2 \rfloor + \; \mbox{deg} \;
a'_{\pi',\mu}\;|\; b_{j, \zeta}^{\pi'}
\neq 0\}. \ee
Hence using (\ref{98o}) in place of (\ref{98}), and then (\ref{119o}) we obtain,
\begin{eqnarray}  \label{scat} 2\mbox{deg} \;\; a_{\pi,\gl}
&\leq& 2\mbox{deg} \;\; a'_{\pi',\mu}+ \cg(\pi') + N - \cg(\gz)\\
&=& 2\mbox{deg} \;\; a'_{\pi',\mu}+ \cg(\pi')  - \cg(\pi) +2.\nonumber
\end{eqnarray}
Therefore by induction
\begin{eqnarray}  \label{scat1} 2\mbox{deg} \;\; a_{\pi,\gl}
&\leq&
2\ell(u)+1 - \cg(\pi) + 2\nonumber\\
&=&
2\ell(w)+1 - \cg(\pi).\nonumber
\end{eqnarray}
giving  condition  (a) in Definition \ref{sog}.
\\ \\
The proof in the case where $\ga$ is an even root is the same as in Section \ref{s.5} apart from the inequalities.  If  $ b_{j,
\zeta}^{\pi'} \neq 0,$ then instead of  (\ref{98}), we have
\begin{eqnarray} \label{98e}
\cg(\zeta) \leq \cg(\pi') + 2(N - j).
\end{eqnarray}
Now condition (a) follows since in place of (\ref{119}) we have, using $m=q=1,$
\begin{eqnarray} \label{1199}\cg(\zeta) = \cg(\pi) + 2(N -1).
\end{eqnarray}
We leave the proof that (b) holds in  Definition \ref{sog} 
to the reader.
\epf 

\section{Powers of \v Sapovalov elements.} \label{zzprod}
\subsection{Isotropic Roots.}\noi First we record  an elementary but important property of the \v Sapovalov element $\theta_{\gamma}$ corresponding to an isotropic root ${\gamma}$.
\bt \label{zprod} If $\gl \in \cH_\gc,$ then $ \theta^2_{\gamma} v_\gl = 0.$ Equivalently,
  $\theta_{\gamma}(\gl-\gc)\theta_{\gamma}(\gl) =0.$\et
\noindent \bpf It is enough to show this for all $\gl$ in a Zariski
dense subset of $\cH_\gc$. We assume $\gc = w\gb$ for $\gb \in \Pi$ where $\gb$ is isotropic, and $w\in  W_{\rm nonisotropic}$.
 Let $\gL$ be the subset of  ${\mathcal H}_{\beta}$ defined by
equation (\ref{tar}), and suppose $\gl \in
w\cdot \Lambda$.
The proof is by induction on the length of $w.$ We can assume that
$w \neq 1.$ Replace $\mu$ with $\mu - \gc'$ and $\gl = s_\ga\cdot \mu$
with $s_\ga\cdot (\mu - \gc') = \gl - \gc$ in Equation (\ref{2nd}) or (\ref{3nd}).
Then $p$ is replaced by $p + q$ and we obtain
\[ e^{p + 2q}_{- \alpha}\theta^{\mu - \gc'}_{\gamma'} = \theta^{\lambda -
\gc}_{\gamma}e^{p + q}_{- \alpha}.\]
 Combining this with Equation
(\ref{2nd}) and using induction we have
\[0 = e^{p+ 2q}_{-\ga}\theta^{\mu - \gc'}_{\gamma'}\theta^{\mu}_{\gamma'} =
\theta^{\gl - \gc}_{\gamma} \theta^\gl_{\gamma}e^{p}_{-\ga} .\] The
result follows since $e_{-\ga}$ is not a zero divisor in
$U({\mathfrak n}^-).$\epf
\noi It follows from the Theorem that for $\gl \in \cH_\gc,$ there is  a sequence of maps
\be
\label{let}  \ldots M(\gl-\gc)
 \stackrel{\psi_{\gl,\gc}}{\longrightarrow} M(\gl)
 \stackrel{\psi_{\gl+\gc,\gc}}{\longrightarrow} M(\gl+\gc) \ldots , \ee
such that the composite of two successive maps is zero.  The map $\psi_{\gl,\gc}$ in (\ref{let}) is defined by $\psi_{\gl,\gc}(x v_{\gl-\gc}) = x\gth_{\gc}v_\gl$.
\subsection{Non-isotropic Roots.}
Now suppose that $\gc$ is non-isotropic. Up to this point we have only evaluated \v{S}apovalov element
$\theta_{\gamma,m}$ at points $\gl \in \cH_{m,\gc}$.
For example under the hypotheses of Lemma \ref{1768} we have,
by \eqref{21a} 
\[ e^{p + mq}_{- \alpha}\theta_{\gamma',m}(\mu)=
\theta_{\gamma,m}(\gl)
e^p_{- \alpha}.\]
In the proof below however, we need to evaluate them at arbitrary points, $\gl \in \fh^*$.  Another situation where it is useful to do this is in the work of Carter \cite{Car} on the construction of orthogonal bases for non-integral Verma modules and simple modules in type A. However if we use the inductive procedure in the proof of Theorems \ref{Shap} and \ref{aShap}, some care must be taken in order that 
evaluation of  $\theta_{\gamma,m}$ at points $\gl$ not in $\cH_{m,\gc}$ is well defined. Indeed if $\fg =\fsl(3)$, and $\gc$ is the positive simple root there are two ways to construct \v{S}apovalov elements using this inductive procedure and their evaluations at such points do not agree.\ff{Recall that \v{S}apovalov elements are only determined modulo a left ideal in $U(\fg)$.}\\ \\
In the proof below, we fix a simple root $\gb$ in the $W$-orbit of $\gc$.  Then starting from $\theta_{\gb,r}(\gl) = e_{-\gb}^r$, if $\gc=w\gb$  we use a fixed reduced expression for $w$ to construct elements 
$\theta_{\gamma,r}(\gl)$ based on \eqref{21a}, 
 for $\gl $ in a dense subset of $\fh^*$.
This determines 
$\theta_{\gamma,r}\in U(\fg)$ which can then be evaluated at arbitrary $\gl\in \fh^*$.

\bt \label{calu} If $\gl \in \cH_{m,\gc}$, $ \theta_{\gamma,m}(\gl)$ is independent of the inductive construction, and  we have
\be\label{CL} 
 \theta_{\gamma,m}(\gl) = 
\theta_{\gamma,1}(\gl-(m-1)\gc)\ldots \theta_{\gamma,1}(\gl-\gc)
\theta_{\gamma,1}(\gl).\ee
\et
\bpf The first statement holds by Theorem \ref{7612}. Clearly \eqref{CL} holds if $\gc$ is a simple root. Suppose that \eqref{21c} holds, and assume that $\gl\in w\cdot \gL$ where $\gL$ is defined in
\eqref{tar}.  For $i=0,\ldots, m-1$ we have $(\mu + \rho-i\gc', \alpha^\vee)= p+iq$, so by the inductive definition 
\be \label{inp} e^{p + (i+1)q}_{- \alpha}\theta_{\gamma',1}(\mu-i\gc')=
\theta_{\gamma,1}(\gl-i\gc)
e^{p+iq}_{- \alpha}.\ee
Combining these equations and using 
the corresponding result for 
$\theta_{\gamma',m}(\mu)$ we obtain
\[  e^{p+mq}_{- \alpha}\theta_{\gamma',m}(\mu) = 
\theta_{\gamma,1}(\gl-(m-1)\gc)\ldots \theta_{\gamma,1}(\gl-\gc)
\theta_{\gamma,1}(\gl)e^{p+iq}_{- \alpha}.\] 
Comparing this to \eqref{21a} and using the fact that $e_{- \alpha}$ is a non-zero divisor, we obtain \eqref{CL}.
 \epf 
\section{An (ortho) symplectic example.} \label{cosp}
\bexa \label{ex2.4}
{\rm  A crucial step in the construction of \v Sapovalov elements was the observation in the proofs of Lemmas \ref{wpfg}  and \ref{wpfg1} that the term $c_{ \zeta,\lambda}$ defined in Equation (\ref{59}) are zero unless $\zeta(\ga) \geq
N-mq,$ (using the notation of the Lemmas).  We give an example where the individual terms on the right of Equation (\ref{59}) are not identically zero, and verify directly that  the sum itself is zero. This cannot happen in Type A.  The key difference in the examples below seems to be that it is necessary to apply Equation (\ref{21nd}) more than once with the same simple root $\ga$. Consider the Dynkin-Kac diagram below for the Lie superalgebra $\fg = \osp(2,4)$.}\eexa

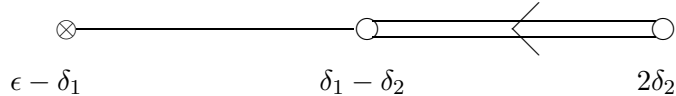
\begin{picture}(310,50)(0,-30)
\thinlines \put(67,0){$\otimes$}
  \put(50,-20){$\epsilon - \delta_{1}$}
  \put(75,3){\line(1,0){105}}
  \put(185,3){\circle{8}}
  \put(167,-20){$\delta_{1} - \delta_{2} $} \put(297,3){\circle{8}}
  \put(287,-20){$2 \delta_{2}$}
  \put(187,0){\line(2,0){107}}
  \put(187,6){\line(2,0){107}}
  \put(240,3){\line(1,-1){10}}
  \put(240,3){\line(1,1){10}}
  \put(335,20){}
\end{picture}

\noi Let $\gb = \epsilon - \gd_1, \ga_1 = \gd_1 - \gd_2, \ga_2 = 2\gd_2$, be the corresponding simple roots.  If we change the grey node to a white node we obtain the Dynkin diagram for $\fsp(6)$.  In this case the simple roots are
$\gb=\delta_{0} - \delta_{1}, \ga_1 = \delta_{1} - \delta_{2} $ and $\ga_2=2 \delta_{2}$. Let
$e_{-\gb}, e_{-\ga_1}, e_{-\ga_2}  $ be the negative simple root vectors.
The computation of  the \v Sapovalov elements
$\gth_1, \gth_2, \gth_3$ for the roots $\gb + \ga_1, \gb + \ga_1+ \ga_2$ and $\gb + 2\ga_1+ \ga_2$ respectively,
is the same for $\osp(2,4)$ and for $\fsp(6)$.
Let $s_1, s_2$ be the reflections corresponding to the simple roots $\ga_1, \ga_2.$ Then
define the other negative root vectors by
$$e_{- \ga_1 - \ga_2} = [e_{-\ga_1},e_{- \ga_2 }  ],  \quad \quad
e_{- 2\ga_1 - \ga_2} = [e_{-\ga_1}, e_{- \ga_1 - \ga_2}],  $$
$$e_{-\gb- \ga_1} =[e_{-\ga_1} ,e_{-\gb}
], \quad e_{-\gb- \ga_1-\ga_2} = [e_{- \ga_2 } ,e_{-\gb- \ga_1}], \quad e_{-\gb- 2\ga_1 - \ga_2} = [e_{-\ga_1},e_{-\gb- \ga_1-\ga_2}]. $$
It follows from the Jacobi identity that
$$[e_{-\gb}, e_{- \ga_1 - \ga_2}] = e_{-\gb- \ga_1-\ga_2},
\quad \quad
[e_{- \ga_1 - \ga_2},e_{-\gb- \ga_1}] = e_{-\gb- 2\ga_1 - \ga_2},$$ and
$$[e_{-\gb},e_{- 2\ga_1 - \ga_2}] = 2e_{-\gb- 2\ga_1-\ga_2}.$$
We order the set of positive roots so that for any  partition $\pi$, $e_{-\ga_1}$ occurs first if at all in $e_{-\pi}$,  and any root vector $e_{-\gs}$ with $\gs$ an odd root occurs last.\\ \\
Suppose that $\gl \in \fh^*$ and define $\gl_1 =s_1\cdot\gl, \gl_2 = s_2\cdot\gl, \mu = s_1\cdot\gl_2$. Let
$$(\gl + \gr,\ga^\vee_1) = p = -(\mu + \gr,(\ga_1+\ga_2)^\vee) $$ and  \[(\gl_1 + \gr,\ga^\vee_2) = (\gl + \gr,(2\ga_1  + \ga_2)^\vee) = q= -(\mu + \gr,(2\ga_1+\ga_2)^\vee) ,\] $$(\gl_2 + \gr,\ga^\vee_1) = (\gl + \gr,(\ga_1+\ga_2)^\vee) = r = -(\mu + \gr,\ga_1^\vee) .$$
Then $r = 2q - p.$ Let $\gc$ be any positive root that involves $\gb$ with non-zero coefficient when expressed as a linear combination of simple roots. We compute the
\v Sapovalov elements
$\gth_{\gc,1}$ for $\fsp(6)$ and $\gth_{\gc}$ for $\osp(2,4).$
To do this we use Equation (\ref{21nd}).   We can assume $\gc \neq \gb.$ Suppose that $p, q, r$ are nonnegative integers. Then
$$e_{-\ga_1}^{p+1}e_{-\gb} = \gth_1e_{-\ga_1}^{p}$$
$$e_{- \ga_2 }^{q+1}\gth_1 = \gth_2e_{- \ga_2 }^{q}$$    $$e_{- \ga_1 }^{r+1}\gth_2 = \gth_3e_{- \ga_1 }^{r}.$$ 
In the computations below we write  $e_{-\pi}$, for $\pi$ a partition (resp. $\gth_i$ for $i=1,2,3$) in place of  $e_{-\pi}v_\gl$ (resp. $\gth_i v_\gl$).
First note that
$$[e_{-\ga_1}^{p+1} ,e_{-\gb}
] = (p+1)e_{-\gb- \ga_1} e^p_{-\ga_1}$$
$$[e_{- \ga_2 }^{q+1} ,e_{-\gb- \ga_1}] =
(q+1)e_{-\gb- \ga_1-\ga_2} e_{- \ga_2 }^q$$
$$[e_{- \ga_2 }^{q+1} ,e_{-\ga_1}]  =-(q+1)e_{- \ga_1 - \ga_2}
e_{- \ga_2 }^q.$$
This easily gives
 $$\gth_1 = (p+1)e_{-\gb- \ga_1} + e_{-\gb}e_{-\ga_1} = pe_{-\gb- \ga_1} + e_{-\ga_1}e_{-\gb}.$$
We order the set of positive roots so that for any  partition $\pi$, $e_{-\ga_2}$ occurs last if at all in $e_{-\pi}$,  and any root vector $e_{-\gs}$ with $\gs$ an odd root occurs first.
\begin{eqnarray}\label{la2}
\gth_2
&=& (p+1)[(q+1)e_{- \gb - \ga_1 - \ga_2} +e_{-\gb- \ga_1}e_{-\ga_2}] +e_{-\gb}[ e_{-\ga_1}e_{- \ga_2} - (q+1)e_{- \ga_1 - \ga_2} ].
\nonumber
\end{eqnarray}
Next order the set of positive roots so that for any  partition $\pi$, $e_{-\ga_1}$ occurs last if at all in $e_{-\pi}$,  and any root vector $e_{-\gs}$ with $\gs$ an odd root occurs first. To find $\gth_3$  we use
$$[e_{-\ga_1}^{r+1} ,e_{-\gb-\ga_1-\ga_2}
] = (r+1)e_{-\gb- 2\ga_1 -\ga_2} e^r_{-\ga_1},$$
$$[e_{-\ga_1}^{r+1} ,e_{-\gb-\ga_1}e_{-\ga_2}
] = (r+1)e_{-\gb- \ga_1}e_{-\ga_1-\ga_2}e^
{r}_{-\ga_1}
+\left( \begin{array}{c}
                r+1 \\
                2 \end{array}\right)e_{-\gb- \ga_1}e_{-2\ga_1 -\ga_2} e^{r-1}_{-\ga_1},$$
$$[e_{-\ga_1}^{r+1} ,e_{-\gb}e_{-\ga_1-\ga_2}
] = (r+1)[e_{-\gb}e_{- 2\ga_1 -\ga_2} e^r_{-\ga_1}
 +e_{-\gb- \ga_1}e_{-\ga_1 -\ga_2} e^r_{-\ga_1} +re_{-\gb-\ga_1}e_{-2\ga_1-\ga_2}
e_{\ga_1}^{r-1}],$$
$$e_{-\ga_1}^{r+1} e_{-\gb}e_{-\ga_2}e_{-\ga_1}
 =
e_{-\gb}[e_{-\ga_2}
e_{- \ga_1}^{2} +
(r+1)e_{- \ga_1-\ga_2}e_{-\ga_1}
+\left( \begin{array}{c}
                r+1 \\
                2 \end{array}\right)e_{- 2\ga_1 -\ga_2}] e^{r}_{-\ga_1}$$
$$+ (r+1)e_{-\gb- \ga_1}[e_{-\ga_2}
e_{- \ga_1} +
re_{- \ga_1-\ga_2}]e^{r}_{-\ga_1}
+(r-1)\left( \begin{array}{c}
                r+1 \\
                2 \end{array}\right)e_{-\gb}e_{- 2\ga_1 -\ga_2} e^{r-1}_{-\ga_1}.$$
The above equations allow us to write $e_{- \ga_1 }^{r+1}\gth_2$  in terms of elements $e_{-\pi}$ with $\pi$ a partition of $\gb + (r+2)\ga_1 + \ga_2$.
We see that the term $e_{-\gb- \ga_1}e_{- 2\ga_1 -\ga_2} e^{r-1}_{-\ga_1}$  occurs in $e_{- \ga_1 }^{r+1}\gth_2$ with coefficient
\[\left( \begin{array}{c}
                r+1 \\
                2 \end{array}\right)[(p+1) -2q + (r-1)] = 0.\]
This is predicted by the cancelation step in the proof of Lemma \ref{wpfg}. In the remaining terms, $e_{- \ga_1 }^{r}$ can be factored on the right, and this yields
\begin{eqnarray}\label{th3}
\gth_3 &=&
(p+1)(q+1)(r+1)e_{-\gb- 2\ga_1 -\ga_2} +(p+1)(q+1)e_{-\gb- \ga_1-\ga_2}e_{-\ga_1} \\
&+& (q+1)(r+1)e_{-\gb- \ga_1}e_{-\ga_1 -\ga_2}-(p/2)(r+1)e_{-\gb}e_{- 2\ga_1 -\ga_2}\nonumber \\
&+& 2(q+1)e_{-\gb- \ga_1}e_{-\ga_2}e_{-\ga_1}+ (r-q+1)e_{-\gb}e_{-\ga_1-\ga_2}e_{-\ga_1}+e_{-\gb}e_{-\ga_2}
e_{- \ga_1}^{2}. \nonumber
\end{eqnarray}
Using the opposite orders on positive roots to those used above to define the $e_{-\pi}$ we obtain
\begin{eqnarray}\label{la3}
\gth_2 &=& p[qe_{- \gb - \ga_1 - \ga_2} +
e_{- \ga_2 }
e_{-\gb- \ga_1}]
+[ e_{-\ga_2}
e_{- \ga_1 } - qe_{- \ga_1 - \ga_2} ]e_{-\gb},
\nonumber
\end{eqnarray}
and
\begin{eqnarray}\label{tho}
\gth_3 &=&
pqr e_{-\gb- 2\ga_1 -\ga_2} +pq
e_{-\ga_1}
e_{-\gb- \ga_1-\ga_2}\\
&+& qr e_{-\ga_1 -\ga_2}e_{-\gb- \ga_1}
-(r/2)(p+1)e_{- 2\ga_1 -\ga_2}e_{-\gb}\nonumber \\
&+& 2q e_{-\ga_1}e_{-\ga_2}e_{-\gb- \ga_1} +
 (r-q-1)e_{-\ga_1}e_{-\ga_1-\ga_2}e_{-\gb}
+e_{- \ga_1}^{2}
e_{-\ga_2}e_{-\gb}.\nonumber
\end{eqnarray}
\br {\rm It seems remarkable that all the coefficients of $\gth_3$ in (\ref{th3}) and (\ref{tho}) are products of linear factors.  This is also true in the Type A case, see Equations (\ref{shtpa}) and (\ref{shtpb}). A partial explanation of this phenomenon is given by specializing these coefficients to zero.
Vanishing of these coefficients gives rise to factorizations of $\gth_3$ as in the examples below.  Factorizations  of \v Sapovalov elements will be discussed in detail elsewhere.
} \er
\bi \itema If $p=(\gl + \gr,\ga_1^\vee) =0$, then $r = 2q$ and we have \begin{eqnarray}\label{la4}
\gth_3 &=&
 2q^2 e_{-\ga_1 -\ga_2}e_{-\gb- \ga_1}
-qe_{- 2\ga_1 -\ga_2}e_{-\gb}\nonumber \\
&+& 2q e_{-\ga_1}e_{-\ga_2}e_{-\gb- \ga_1} +
 (q-1)e_{-\ga_1}e_{-\ga_1-\ga_2}e_{-\gb}
+e_{- \ga_1}^{2}
e_{-\ga_2}e_{-\gb} \nonumber \\&=& \gth_{\ga_1+\ga_2}\gth_{\gb+\ga_1}.\nonumber
\end{eqnarray}
\itemb
If $q=(\gl + \gr,(2\ga_1  + \ga_2)^\vee) =0$ then $p = -r$, $\gth_2=
\gth_{\gb +\ga_1+\ga_2} = \gth_{\ga_2}\gth_{\gb+\ga_1}$, and we have \begin{eqnarray}\label{la5}
\gth_3 &=&
[(p/2)(p+1)e_{- 2\ga_1 -\ga_2}
 -(p+1)e_{-\ga_1}e_{-\ga_1-\ga_2}
+e_{- \ga_1}^{2}
e_{-\ga_2}]e_{-\gb}\nonumber \\&=& \gth_{2\ga_1+\ga_2}\gth_{\gb}.\nonumber
\end{eqnarray}
\itemc If $r=(\gl + \gr,(\ga_1  + \ga_2)^\vee) = 0$, then $p = 2q$,
 and
we have \begin{eqnarray}\label{la6}
\gth_3 &=&
2q^2
e_{-\ga_1}
e_{-\gb- \ga_1-\ga_2}\nonumber \\
&+& 2q e_{-\ga_1}e_{-\ga_2}e_{-\gb- \ga_1} -
 (q+1)e_{-\ga_1}e_{-\ga_1-\ga_2}e_{-\gb}
+e_{- \ga_1}^{2}
e_{-\ga_2}e_{-\gb}\nonumber \\
 &=&
e_{-\ga_1}[2q^2
e_{-\gb- \ga_1-\ga_2}
 2q e_{-\ga_2}e_{-\gb- \ga_1} -
 (q+1)e_{-\ga_1-\ga_2}e_{-\gb}
+e_{- \ga_1}
e_{-\ga_2}e_{-\gb}]\nonumber \\&=& \gth_{\ga_1}\gth_{\gb+\ga_1+\ga_2}.\nonumber
\end{eqnarray}
\ei
Similarly if $p=-1,$ (resp. $q=-1,$ $r=-1$) then (\ref{tho}) yields the factorizations $\gth_3=\gth_{\gb+\ga_1}\gth_{\ga_1+\ga_2}$, (resp. $\gth_2=
 \gth_{\gb+\ga_1}
\gth_{\ga_2},$ $\gth_3=\gth_{\gb}\gth_{2\ga_1+\ga_2}$, and  $\gth_3=\gth_{\gb+\ga_1+\ga_2}\gth_{\ga_1}$). On the other hand we see that $p$ divides the coefficients of $e_{-\gb- 2\ga_1 -\ga_2}$ and $e_{-\ga_1}
e_{-\gb- \ga_1-\ga_2}$ in (\ref{tho}) since when $p = 0$, $\gth_3=\gth_{\ga_1+\ga_2}\gth_{\gb+\ga_1}$ can be written as a linear combination of different $e_{-\pi}$. In this way we obtain explanations for all the linear factors in (\ref{th3}) and (\ref{tho}) with the exception of the coefficients $r-q\pm 1$ of
$e_{-\gb}e_{-\ga_1-\ga_2}e_{-\ga_1}$ and $e_{-\ga_1}e_{-\ga_1-\ga_2}e_{-\gb}.$
At this point it may be worthwhile mentioning that $r-q =(\gl + \gr,\ga_2^\vee)$.  In addition equality holds in the upper bounds given in Theorem \ref{Shap} for the degrees of all the coefficients in (\ref{th3}) and (\ref{tho}).
\section{The Type A Case.}\label{s.8}
\subsection{Lie Superalgebras.}\label{far} We construct the elements $\theta_{\gamma}$ in Theorem \ref{Shap}
explicitly when ${\mathfrak g} = \fgl(m,n).$  Suppose that $\gc = \gep_r
- \gd_{s}.$ For $1 \leq  i<j \leq  m$ and $1 \leq  k<\ell \leq  n$
define roots $\gs_{i,j}, \gt_{k,\ell}$ by
\[ \gs_{i,j} = \gep_{i} - \gep_{j}, \quad \gt_{k,\ell} = \gd_k - \gd_\ell. \]
Suppose   $B = (b_{i,j})$ is a $k \times \ell$ matrix with entries
in $U({\mathfrak n}^{-})$ , $I \subseteq \{1, \ldots , k \}$ and $J
\subseteq \{1, \ldots , \ell \}.$  We denote the submatrix of $B$
obtained by deleting the $ith$ row for $i \in I,$ and the $jth$
column for $j \in J$ by $_IB_J.$ If either set is empty, we omit the
corresponding subscript. When $I = \{ i\},$ we write $_iB$ in place
of $_IB$  and likewise when $|J| = 1.$

If   $k = \ell$ we define two noncommutative determinants of $B,$
the first working from left to right, and the second working from
right to left.
$${\stackrel{\longrightarrow}{{\rm det}}}(B) =  \sum_{w \in \mathcal{S}_k} sign(w) b _{w(1),1} \ldots b _{w(k),k} ,$$ $$ {\stackrel{\longleftarrow}{{\rm det}}}(B) =  \sum_{w \in \mathcal{S}_k} sign(w) b _{w(k),k} \ldots
b _{w(1),1} .$$ If $k = 0$ we make the convention that
${\stackrel{\longleftarrow}{{\rm det}}}(B) =
{\stackrel{\longrightarrow}{{\rm det}}}(B) = 1.$ Consider the
following matrices with entries in $U(\fn^-)$
  $$
A^+(\gl, r) =  \left[ {\begin{array}{ccccc}
e_{r+1,r} &  e_{r+2,r} & \hdots & e_{m,r}\\
-a_1 & e_{r+2,r+1} & \hdots & e_{m,r+1} \\
0&-a_2 & \hdots & e_{m,r+2} \\
 \vdots  & \vdots & \ddots & \vdots              \\
0&  \hdots & \ldots&
e_{m,m-1}   \\
0 & \hdots & 0& -a_{m-r}\\
\end{array}}
 \right],$$

$$A^-(\gl, s) =  \left[ {\begin{array}{cccc}
b_1& 0 & \hdots & 0   \\
e_{m+2,m+1} & b_2& \hdots & 0   \\
 \vdots  & \vdots & \ddots & \vdots             \\
 e_{m+s-1,m+1} & e_{m+s-1,m+2} & \hdots & b_{s-1}\\
e_{m+s,m+1} & e_{m+s,m+2} & \hdots & e_{m+s,m+s-1}    \\
\end{array}}
 \right],
$$ where $a_j=(\gl + \gr,\gs_{r,r+j}^\vee ) $ and $b_i = (\gl + \gr,\gt_{i,s}^\vee )$.
Also  let$B^+(\gl, r)$ (resp. $B^-(\gl, s)$) be the matrices obtained from $A^+(\gl, r)$ (resp. $A^-(\gl, s)$) by replacing each $a_i$ by $a_i -1$ (resp. replacing each $b_i$ by $b_i -1$).
Observe that in $A^+(\gl, r)$ and $A^-(\gl, s)$ the number of rows exceeds the number of columns by one. We also consider two degenerate cases:
if $r=m$ then $A^+(\gl, r)$
and $B^+(\gl, r)$ are ``matrices with zero rows"
  In this case we ignore the summation over $j$ in the following formulas,  replacing
${\stackrel{\longrightarrow}{{\rm det}}}(_jA^+(\gl, r))$ by 1 and $i+j+r+m$ by $i+1+r+m$. Similar remarks apply to the case where $s = 1.$
 \\ \\
Below we present a
two determinantal formulas
for the \v Sapovalov element $\gth_\gc$ evaluated at $\gl\in \cH_\gc$.  A key difference is the placement of the odd root vectors $e_{m+i,j+r-1}$.
 \bt \label{shgl}
\be \label{shtpa} \theta_{\gamma}(\gl)
= (-1)^{m+r}{\sum_{j=1}^{m-r+1} \sum_{i= 1}^s
(-1)^{i+j}{\stackrel{\longrightarrow}{{\rm det}}}(_jA^+(\gl, r))
\;{\stackrel{\longleftarrow}{{\rm det}}}( _iA^-(\gl, s))}e_{m+i,j+r-1}.\quad
\quad \ee
\be \label{shtpb}= {\sum_{j=1}^{m-r+1} \sum_{i= 1}^s e_{m+i,j+r-1}\;(-1)^{i+j+r+m}{\stackrel{\longleftarrow}{{\rm det}}}(_jB^+(\gl, r))\;{\stackrel{\lra}{{\rm det}}}( _iB^-(\gl, s))}.\quad \quad \ee

\et
\bpf First we prove (\ref{shtpa}).
 For the isotropic simple root $\gb = \gep_m-\gd_1$, (\ref{shtpa}) reduces to $\gth_\gb(\gl)  = e_{m+1, m}.$ 
Note that the overall sign in \eqref{shtpa} is determined by the condition that the coefficient of $e_{-\pi^0}$ in $\theta_{\gamma}(\gl)$ is equal to 1, and that this term arises when the first row of  $A^-(\gl, s)$ and the last row of 
$A^+(\gl, r))$ are deleted before taking determinants. 
Suppose that $\ga = \gd_{s} -\gd_{s+1}, \gc =
\gep_r - \gd_{s}$ and $\gc' = \gep_r - \gd_{s+1} = s_\ga \gc.$ 
Assuming the result for $\gc$ we prove it for $\gc'$.  The result
for $\gep_{r -1} - \gd_{s}$ can be deduced in a similar way. Set
$e_{-\ga}
= e_{m+s+1,m+s}.$ For $1 \leq i \leq s-1$ we have 
\[ \gt_{i,s+1} =  \gd_i - \gd_{s+1} = s_\ga \gt_{i,s}.\]  
Consider the matrix
$$A^-(\gl',s+1) =  \left[ {\begin{array}{cccc}
-(\gl' + \gr,\gt_{1,s+1}^\vee )  & 0 & \hdots & 0   \\
e_{m+2,m+1} & -(\gl' + \gr,\gt_{2,s+1}^\vee )  & \hdots & 0   \\
 \vdots  & \ddots & \vdots      &\vdots       \\
 e_{m+s-1,m+1} & e_{m+s-1,m+2} & \hdots  & 0   \\
 e_{m+s,m+1}  &  e_{m+s,m+2}& \ldots &-(\gl' + \gr,\gt_{s,s+1}^\vee ) \\
e_{m+s+1,m+1} & e_{m+s+1,m+2} & \hdots &e_{m+s+1,m+s}    \\
\end{array}}
 \right].
$$
The matrix $A^-(\gl',s+1)$  replaces the matrix $A^-(\gl, s)$ in
Equation (\ref{shtpa}) in the analogous expression for
$\theta_{\gc'}(\gl')$. Suppose that $(\gl + \gr,\ga^\vee )= p$ and
let $\gl' = s_\ga\cdot\gl.$  Then
\[ (\gl' + \gr,\gt_{i,s+1}^\vee ) = (\gl + \gr,\gt_{i,s}^\vee ),\] for $1 \leq i \leq s-1$
, and this means that the first $s-1$ diagonal entries of
$A^-(\gl',s+1)$ and $A^-(\gl, s)$ are equal.  Also
the entry in row $s$ and column $s$ of
$A^-(\gl',s+1)$ is equal to $-p.$
If we remove the last row (row $s+1$) from $A^-(\gl',s+1)$  the last column of the resulting matrix will have only one non-zero entry $-p$.
If in addition we remove this column, we obtain the matrix $A^-(\gl, s)$. Therefore
 \be \label{stra2} -p\;{\stackrel{\longleftarrow}{{\rm det}}}(
{}_sA^{-}(\gl,s)) = {{\stackrel{\longleftarrow}{{\rm det}}}}( {}_{s+1} A^-(\gl',s+1)).\ee
Similarly by removing row $s$ from $A^-(\gl',s+1)$
and noting that $e_{m+s+1,m+s} $ commutes with all entries in $
{}_sA^{-}(\gl,s)$,
we see that
\be \label{stra1}
e_{m+s+1,m+s} {\stackrel{\longleftarrow}{{\rm det}}}(
{}_sA^{-}(\gl,s))=
{\stackrel{\longleftarrow}{{\rm det}}}(
{}_sA^{-}(\gl,s)) e_{m+s+1,m+s} \quad  =
{\stackrel{\longleftarrow}{{\rm det}}}( _{s} A^-(\gl',s+1)). \ee
Equation (\ref{21a})  in this situation takes the form
 \begin{equation} \label{22nd}
e^{p + 1}_{- \alpha}\theta_\gc(\gl) = \theta_{\gc'}(\gl') e^p_{-
\alpha}.
\end{equation}
If   $r \leq  k \leq m +s -1$ we have
 \begin{equation} \label{23nd}
e^{p + 1}_{- \alpha} e_{m+s,k} = (pe_{m+s+1,k} +
 e_{m+s+1,m+s} e_{m+s,k})e^p_{- \alpha}.
\end{equation}
We now consider  two cases: in the first entries
in ${\stackrel{\longleftarrow}{{\rm
det}}}( _iA^-(\gl, s))$ are
replaced by entries in
${\stackrel{\longleftarrow}{{\rm det}}}(
_iA^-(\gl',s+1))$.
Suppose $1 \leq  i \leq  s-1$, and $r \leq  k \leq m$.
 Then $
e_{- \alpha}$ commutes with $e_{m+i,k}$ and all entries in the
matrix $_iA^-(\gl, s)$ except for those in the last row. Replacing
$e_{m+s,k}$ in the matrix
$_iA^-(\gl, s)$ by $e_{m+s+1,k}$
yields the matrix $_{\{i,s\}}A^-(\gl',s+1)_s.$ Hence Equation
(\ref{23nd}) gives the first equality below.  For the second we
use a cofactor expansion,
\begin{eqnarray} \label{24nd}
e^{p + 1}_{- \alpha}\;{\stackrel{\longleftarrow}{{\rm
det}}}( _iA^-(\gl, s))e_{m+i,k} &=&
[p\;{\stackrel{\longleftarrow}{{\rm det}}}(
_{\{i,s\}}A^-(\gl',s+1)_s)
\nonumber \\& +& e_{m+s+1,m+s}
 {\stackrel{\longleftarrow}{{\rm det}}}(
_iA^-(\gl, s)) ]e_{m+i,k}
e^p_{- \alpha}\nonumber\\
&=& \;{\stackrel{\longleftarrow}{{\rm det}}}(
_iA^-(\gl',s+1))e_{m+i,k} e^p_{- \alpha}.
\end{eqnarray}
Now consider the case $i=s.$ Here
entries in
${\stackrel{\longleftarrow}{{\rm
det}}}( _sA^-(\gl, s))$
 are unchanged but the factor $e_{m+s,k}$ is replaced.
If $r \leq  k \leq m$, then all
entries in the matrix $_sA^-(\gl, s),$ commute with $e_{m+s,k}$ and
$e_{- \alpha}, $ so by Equation (\ref{23nd}) we get the first
equality below, and the  second equality comes from Equations
(\ref{stra2}) and (\ref{stra1})
\be\label{2576}
e^{p + 1}_{- \alpha} \;{\stackrel{\longleftarrow}{{\rm
det}}}( _sA^-(\gl, s))e_{m+s,k}e^{-p}_{- \alpha}=
\;{\stackrel{\longleftarrow}{{\rm det}}}( {}_sA^{-}(\gl,s))
[pe_{m+s+1,k} +
 e_{m+s+1,m+s}e_{m+s,k}] \nonumber \ee
\be\quad\quad\quad\quad\quad\quad\quad\quad\quad = -{\stackrel{\longleftarrow}{{\rm det}}}( _{s+1}A^-(\gl',s+1))e_{m+s+1,k}  + {\stackrel{\longleftarrow}{{\rm det}}}(_{s} A^-(\gl',s+1))e_{m+s,k}.\ee
 Since ${\stackrel{\longrightarrow}{{\rm
det}}}(_jA^+(\gl, r))$ commutes with $e_{- \alpha}$ and $e_{m+j,i}$
for all $i,j$, it follows from the induction assumption and Equations (\ref{22nd}),
(\ref{24nd}) and (\ref{2576}) that \be \label{shtpab}
\theta_{\gc'}(\gl')= (-1)^{m+r}\sum_{j=1}^{m-r+1} \sum_{i= 1}^{s+1}
 (-1)^{i+j}{\stackrel{\longrightarrow}{{\rm
det}}}(_jA^+(\gl, r)) \;{\stackrel{\longleftarrow}{{\rm det}}}(
_iA^-(\gl',s+1))e_{m+i,j+r-1}.\quad \quad \nonumber \ee as desired.  The proof of  (\ref{shtpb}) is similar using \eqref{21b} in place of \eqref{21a}.\epf

\subsection{A variation.}\label{sov}

By inserting odd root vectors as an extra column in either
$A^+(\gl, r)$  or $A^-(\gl, s)$, we obtain a variation of Equations \eqref{shtpa} 
 where the odd root vectors are inserted into one of the determinants. Note that the two determinants in \eqref{shtpa} commute.  We give the details only for $A^+(\gl, r)$.
For $1\le i\le s$, let $C^{(i)}(\gl, r)$ be the matrix obtained from $A^+(\gl, r)$ by adjoining the vector
\[ (e_{m+i,r}, e_{m+i,r+1}, \ldots, e_{m+i,m} )^\trp\]
as the last column.

 \bt \label{thgl}With the above notation we have
 \be \label{thtpa} \theta_{\gamma}(\gl)
=
\sum_{i= 1}^s
(-1)^{i+1}{\stackrel{\longleftarrow}{{\rm det}}}( _iA^-(\gl, s))\;{\stackrel{\longrightarrow}{{\rm det}}}(C^{(i)}(\gl, r))
\quad
\quad \ee

\et
\bpf This follows by cofactor expansion of
${\stackrel{\longrightarrow}{{\rm det}}}(C^{(i)}(\gl, r))$ down the last column.
\epf

\subsection{Lie Algebras.}
Let $\fg= \fgl(m)$, and $\ga = \gep_r -\gep_t$.  We give a determinantal formula for the  \v{S}apovalov element $\gth_{\ga,1}$. Set $\gs_{i,j} = \gep_{i} - \gep_{j}$. Consider the
following matrix with entries in $U(\fn^-)$.

 $$
C_\gl=  \left[ {\begin{array}{cccccc}
e_{r+1,r} &  e_{r+2,r} & \hdots & e_{t-1,r}& e_{t,r}\\
-(\gl + \gr,\gs_{r,r+1}^\vee ) & e_{r+2,r+1} & \hdots &e_{t-1,r+1}& e_{t,r+1} \\
0&-(\gl + \gr,\gs_{r,r+2}^\vee ) & \hdots &e_{t-1,r+2}& e_{t,r+2} \\
 \vdots  & \vdots & \ddots & \vdots &\vdots              \\
0&  \hdots & \ldots&
-(\gl + \gr,\gs_{r,t-1}^\vee )&e_{t,t-1}   
\end{array}}
 \right].$$

\bt \label{shgl2}
The \v{S}apovalov element $\gth_{\ga,1}$ is given by
\[ \gth_{\ga,1}(\gl) =
{\stackrel{\longrightarrow}{{\rm
det}}}(C_\gl), \]
for $\gl$ such that $ (\gl + \gr,\gs_{r,t}^\vee )=1$.\et
\bpf Similar to the proof of Theorem \ref{shgl}. \epf
\noi Note that if $ (\gl + \gr,\gs_{r,t}^\vee )=1$ and $ (\gl + \gr,\gs_{r,s}^\vee )=0$, then $ (\gl + \gr,\gs_{s,t}^\vee )=1$ and $ (\gl + \gr-
\gs_{s,t},\gs_{r,s}^\vee )=1$.
\bc \label{sapc}  In the above situation, for any highest weight vector $v_\gl$ of weight $\gl$, we have
$\gth_{\gep_{r,1} -\gep_{t,1}}v_\gl = 
\gth_{\gep_{r,1} -\gep_{s,1}} \gth_{\gep_{s,1} -\gep_{t,1}}v_\gl$.
\ec
\bpf Under the given hypothesis the matrix $C_\gl$ in Theorem \ref{shgl2}  is  block upper triangular.\epf
\noi
\bc \label{1.5} For $p\ge 1$, we have
\[\gth_{\gc,p}(\gl) = C_{\gl-(p-1)\gc}\ldots C_{\gl-\gc}C_{\gl}.\]
\ec
\bpf Combine Theorems \ref{calu} and \ref{shgl2}. \epf
\noi  Determinants similar to those in Theorem \ref{shgl2} are introduced in \cite{CL} Equation (5).  Corollary  \ref{1.5} may be viewed as a version of \cite{CL} Theorem 2.7. Note however that Carter and Lusztig work inside the Kostant $\Z$-form of $U(\fg).$ They use their result to study tensor powers of the defining representation of $GL(V)$, and homomorphisms between Weyl modules in positive characteristic.   See also \cite{Br3}, \cite{Carlin}, \cite{Car} and \cite{F}.


\section{{Survival of \v{S}}apovalov elements in factor modules.}\label{surv}
 Let  $v_{\gl}$ be a highest weight vector
in a Verma module $M(\gl)$ with highest weight $\gl,$ and suppose $\gc$ is an odd root with
$(\gl+\gr, \gc) =0$. We are interested in the condition that the image of
$\theta_{\gc}v_{\gl}$ is non-zero in various factor modules of $M(\gl)$.
\subsection{Independence of \v{S}apovalov elements.}
Given $\gl \in \fh^*$ recall the set $B(\gl)$ defined in Section \ref{s.1}, and define a ``Bruhat order" $\le$ on $B(\gl)$ by $\gc' \le \gc$ if $\gc-\gc'$ is a sum of positive even roots. Then introduce a relation $\downarrow$ on $B(\gl)$  by $\gc' \da \gc$ if $\gc' \le \gc$ and $(\gc,\gc') \neq 0.$ If $\gc \in B(\gl)$, we say that $\gc$ is $\gl$-{\it minimal} if $\gc' \da \gc$  with $\gc' \in B(\gl)$ implies that $\gc' = \gc$.
For $\gc \in B(\gl)$ set $B(\gl)^{-\gc} = B(\gl)\backslash \{\gc\}$. We say $\gc$ is {\it independent
at } $\gl$ if $$\theta_{\gc}v_{\gl}\notin \sum_{\gc' \in B(\gl)^{-\gc} }U(\fg)\theta_{\gc'} v_{\gl}.$$
\bp \label{pot}If $\gc' \da \gc$ with $\gc' \in B(\gl)$ and $\gc' < \gc$, then $\theta_{\gc}v_{\gl}\in U(\fg)\theta_{\gc'} v_{\gl}.$ \ep
\bpf The hypothesis implies that
$(\gc,\ga^\vee)>0$
and $\gc =s_\ga \gc'$. Thus the   result follows from \cite{M17}.\epf
\noi By the  Proposition,  if we are interested in the independence of the \v{S}apovalov elements $\gth_\gc$ for distinct isotropic roots, it suffices to  study only $\gl$-minimal roots $\gc$.
\\ \\
For the rest of this section we assume that $\fg = \fgl(m,n)$. We use 
Equation
(\ref{shtpa}), and  order the positive roots of $\fg$ so that each summand in this equation is a constant multiple of $e_{-\pi}$ for some $\pi \in {\overline{\bf P}}(\gamma)$. For such $\pi$  the odd root vector is the rightmost factor of $e_{-\pi}$, that is we have $e_{-\pi}\in U(\fn^-_0)\fn^-_1$.
\bl \label{car}If If $\gc$ is $\gl$-minimal, then $e_{-\gc} v_\gl$ occurs with non-zero coefficient in $\gth_\gc v_\gl$.
\el
\bpf Assume $\gc = \gep_r
- \gd_{s}.$ Then if $\ga = \gep_r - \gep_i$ with $r<i$,  or $\ga = \gd_j - \gd_s$ with $j<s$ we have $(\gl + \gr, \ga^\vee) \neq 0$, since $\gc$ is $\gl$-minimal.  In other words the entries on the superdiagonals of $A^+(\gl, r)$ and $A^-(\gl, s)$ are non-zero.  Thus the result follows from Theorem \ref{shgl}.\epf
\bt \label{boy} The isotropic root $\gc$ is independent at  $\gl$ if and only if $\gc$ is $\gl$-{minimal}. \et
\bpf Set $B=B(\gl)^{-\gc} $. If
$\gc$ is not $\gl$-{minimal}
then $\gc$ is not independent at  $\gl$
by Proposition \ref{pot}. Suppose that $\gc$ is $\gl$-{minimal}
and
$$\theta_{\gc}v_{\gl}\in \sum_{\gc' \in B}U(\fg)\theta_{\gc'} v_{\gl} = \sum_{\gc' \in B}U(\fn^-)\theta_{\gc'} v_{\gl},$$ then by comparing weights
\be \label{la}\theta_{\gc}v_{\gl}\in \sum_{\gc' \in B}U(\fn^-_0)e_{-\gc'} v_{\gl},\ee
But
Lemma \ref{car} implies that
$$\theta_{\gc}v_{\gl}\equiv ce_{-\gc}v_{\gl}\mod \sum_{\gc' \in B}U(\fn^-_0)e_{-\gc'} v_{\gl}$$ for some non-zero constant $c.$  Combined with
(\ref{la}) we obtain a contradiction to the PBW Theorem.
\epf
\subsection{Survival of \v{S}apovalov elements in Kac modules.}
For $\fg = \fgl(m,n)$ we have $\mathfrak{g}_1 = \mathfrak{g}_1^+ \oplus
\mathfrak{g}_1^-,$ where $\mathfrak{g}_1^+ $ (resp. $\mathfrak{g}_1^-$) is the set of block  upper (resp. lower) triangular matrices.
Let $\fh$ be the Cartan subalgebra of $\fg$ consisting of diagonal matrices, and set $\mathfrak{p} =  \mathfrak{g}_0  \oplus \mathfrak{g}_1^{+}.$ Next let
 \begin{eqnarray*}
  P^+ &=&
   \{ \lambda \in \mathfrak{h}^*|(\lambda, \alpha^\vee) \in \mathbb{Z}, (\lambda, \alpha^\vee) \geq 0 \quad
 \mbox{for all}\quad \alpha \in \Delta^+_0 \}
 \end{eqnarray*}
For $\gl \in P^+,$ let $L^0(\lambda)$ be the (finite dimensional) simple  $\mathfrak{g}_0$-module with highest weight $\gl$.  Then
$L^0(\lambda)$ is naturally a $\fp$-module and we define the Kac module $K(\gl)$ by
\[K(\lambda) =  U(\mathfrak{g}) \otimes_{U(\fp)}
L^0(\lambda).\]
 Note that as a $\mathfrak{g}_0$-module
\[K(\lambda) =  \Lambda(\mathfrak{g}^-_1) \otimes
L^0(\lambda).\]
The next result is well-known.  Indeed two
methods of  proof are given in  Theorem 4.37 of \cite{Br}.  The
second of these is based on Theorem 5.5 in \cite{S2}.  We give a short
proof using  Theorem \ref{shgl}. We now assume that the roots are ordered as in Equation (\ref{shtpb}), that is  with the odd root vector first.
\bt \label{shapel} If $\lambda$ and $\lambda - \epsilon_r+\delta_s$
belong to $P^+$ and $(\lambda+\rho,\epsilon_r - \delta_s) = 0$, then
$$[K(\lambda):L(\lambda-\epsilon_r+\delta_s)] \neq 0.$$ \et
 \bpf Set $\gamma = \epsilon_r -
\delta_s$.
 Let $\theta_\gamma(\lambda)$ be as in Theorem \ref{shgl}.  Then
$w = \theta_\gamma(\lambda)v_\lambda$ is a highest weight vector in
the Verma module $M(\lambda)$ with weight $\lambda - \gamma$.  It
suffices to show that the image of $w$ in the Kac module
$K(\lambda)$ is nonzero.  We have an embedding of $\fg_0$-modules
\[ \fg^-_1 \otimes L^0(\lambda) \subseteq \Lambda \fg^-_1 \otimes L^0(\lambda). \]
The elements $e_{m+j,i+r-1}$ in Equation (\ref{shtpb}) form part of
a basis for $\fg^-_1$. Furthermore the coefficient of
$e_{m+j,i+r-1}$ belong to $U(\fn^-_0)$. Therefore it suffices to
show that the coefficient of $e_{m+s,r}$ in this equation is
nonzero. This coefficient is found by deleting the first column  of the matrix $B^+(\gl, r)$ and the last row of
$B^-(\gl, s)$ and taking determinants of the resulting matrices, which have only zero entries above the main diagonal. We find that the coefficient of $e_{m+s,r}$ is
\[  \pm \prod_{k=1}^{m-r} (1 - (\lambda + \rho, \gs_{r,r+k}^\vee))
\prod^{s-1}_{k=1} (1 - (\lambda + \rho, \gt_{k,s}^\vee ))    \]
Since $\lambda \in P^+, (\lambda + \rho, \gs_{r,r+k}^\vee) \geq 1$
with equality if and only if $k = 1$ and $(\lambda, \epsilon_r -
\epsilon_{r+1}) = 0$.    This cannot happen if $\lambda - \gamma \in
P^+$, so the first product above is nonzero, and similarly so is the
second. \epf
\section{Changing the Borel subalgebra.}\label{sscbs}
\subsection{Adjacent Borel subalgebras.}
 We consider the behavior of \v Sapovalov elements when the Borel subalgebra is changed.
Let ${\fb'},{\fb}''$ be arbitrary adjacent Borel subalgebras, 
 and suppose
\be \label{fggb} \fg^\ga \subset {\fb}', \quad \fg^{-\ga} \subset {\fb''}\ee for some isotropic root $\ga$.
Let $S$ be the intersection of the sets of roots of  $\fb'$ and $\fb''$, $\fp =\fb' + \fb''$ and $\fr  = \bigoplus_{\gb \in S} \ttk e_{-\gb}.$ Then $\fr, \fp$ are subalgebras of $\fg$ with $\fg = \fp \oplus \fr$.  Furthermore $\fr$ is stable under $\ad e_{\pm \ga},$ and consequently, so is $U(\fr)$.
 Note that  \be \label{rgb}\rho({\fb}'') = \gr({\fb'}) + \ga.\ee
Also if $v_\mu$ is a highest weight vector with weight  $\mu$, then \be \label{ebl} e_\ga e_{-\ga} v_\mu = h_\ga v_\mu =(\mu, \ga)v_\mu . \ee
As is well known (see for example Corollary 8.6.3 in \cite{M}),  if $(\mu +\rho({\fb'}), \ga) \neq 0$, then
\be \label{mon}\mu({\fb''}) + \rho({\fb''}) = \mu({\fb}') + \rho({\fb}').\ee
In this situation we call the change of Borel subalgebras from ${\fb'}$ to ${\fb}''$ (or vice-versa) a {\it typical change of Borels}.
Consider the Zariski dense subset $\gL_{\gc, m}$ of ${\mathcal H}_{\gc, m}$ given by
\[\gL_{\gc, m} =  \{ \mu \in
{\mathcal H}_{\gc, m}|(\mu + \rho, \ga) \notin \bbZ  \mbox{ for all positive roots } \ga \neq \gc\}.\]
When $\gc$ is isotropic, we write $\gL_{\gc} $ in place of $\gL_{\gc, 1} $. Since the coefficients of
$\theta_{\gamma,m}$ are polynomials,
$\theta_{\gamma,m}$ is determined (as usual modulo a left ideal) by the values of
$\theta_{\gamma,m}v_\mu$ for
$\mu \in \gL_{\gc, m}$.
\\ \\
Suppose that $\gc$ is a positive root of both ${\fb'}$ and ${\fb}''$, and that
$\theta_{\gamma,m}$ is a \v Sapovalov element
 corresponding to the pair $(\gc, m)$ using the negatives of the roots of ${\fb}''$. For brevity set $\gth =\theta_{\gamma,m}(\mu).$ Assume that $v_{\mu'}$ is a highest weight  vector in a Verma  module $M_{\sfb'}(\mu')$ for ${\fb'}$ with highest weight  ${\mu'}\in \gL_{\gc, m}.$
Then $v_\mu =e_{-\ga} v_{\mu'}$ is  a highest weight  vector for ${\fb}''$
which also generates $M_{{ \sfb'}}(\mu')$.  Thus we can write
\[M_{{ \sfb'}}(\mu')=M_{{ \sfb''}}(\mu).\]
Next note that $u = \theta e_{-\ga} v_{\mu'}$ is  a highest weight  vector for ${\fb}''$, and $e_{\ga} \theta e_{-\ga} v_{\mu'} $ is a highest weight vector for ${\fb'}$ of weight $\mu' -m\gamma$ that generates the same submodule of $M_{{ \sfb'}}(\mu')$ as $u$.
We can write $\gth$ in a unique way as $\theta=e_\ga \theta_{1} + \theta_{2}$  with $\gth_i \in U(\fr)$. Then
\begin{eqnarray} \label{slab}
 e_{\ga}\theta e_{-\ga} v_{\mu'} &=& e_{\ga}\theta_2 e_{-\ga} v_{\mu'} \nonumber\\
 &=& \gth'_1 e_{-\ga} v_{\mu'} \pm  \theta'_2 v_{\mu'} \nonumber
 \end{eqnarray}
where $\gth'_1 = [e_{\ga},\theta_2], \; \; \gth'_2 = (\mu',\ga)\theta_2 \in U(\fr).$
 Note that the term $e_{-m\pi_\gc}$ cannot occur in $e_\ga \gth_1$ or $\gth'_1 e_{-\ga}. $ 
  Allowing for possible re-ordering of positive roots used to define the $e_{-\pi}$ (compare Lemma \ref{delorder}) we conclude that modulo terms of lower degree, the coefficient of $e_{-m\pi_\gc}$ in
 $e_{\ga}\theta e_{-\ga} v_{\mu'}$ is equal to
 $\pm(\mu',\ga)$
 times the coefficient of
 $e_{-m\pi_\gc}$
  in  $\theta_2 e_\ga v_{\mu'}$. Since ${\mu'}\in \gL_{\gc, m},$ each change of Borels in (\ref{distm}) is typical, and thus the foregoing applies to each link in the chain.
\subsection{Chains of Borel subalgebras.}
Using adjacent Borel subalgebras it is possible to give an alternative construction of \v Sapovalov elements corresponding to an isotropic root $\gc$ which is a  simple root for some Borel subalgebra.  This condition always holds in type A, but for other types,  it is quite restrictive:  if $\fg = \osp(2m,2n+1)$ the assumption only holds
for roots of the form $\pm(\gep_i-\gd_j)$, while if $\fg = \osp(2m,2n)$ it holds only
for these roots and the root $\gep_m+\gd_n$. (Theorem \ref{aShap} on the other hand applies to any positive isotropic root, provided we choose the appropriate Borel subalgebra satisfying Hypothesis (\ref{i}).)
\\ \\
Suppose that  
$\fb$ is the distinguished or anti-distinguished  Borel subalgebra, and let $\mathfrak{b}'$
be another Borel subalgebra with the same even part as
$\mathfrak{b}.$ Consider  a  sequence of Borel subalgebras
\[\mathfrak{b} = \mathfrak{b}^{(0)}, \mathfrak{b}^{(1)}, \ldots,
 \mathfrak{b}^{(r)} \] as in Equation (\ref{distm}).
  There are isotropic roots $\ga_i$ such that $\fg^{\ga_i} \subset \fb^{(i-1)}, \fg^{-\ga_i} \subset \fb^{(i)}$  for $1 \leq i \leq r$, and  $\ga_1,\ldots,\ga_r$ are distinct positive roots of $\fb.$
Set $v_0 = v_\gl, \gl_{0} = \gl$, and then $v_i= e_{-\ga_i}v_{i-1}, \gl_{i} = \gl_{i-1}-\ga_i$ for $1 \leq i \leq r$. Next set $u_r =e_{-\gc}v_r$ and 
$u_i = e_{\ga_{i+1}}\ldots e_{\ga_r}u_r$, so that $u_i$ and $v_i$ are  highest weight vectors for the Borel subalgebra $\fb^{(i)}$ with weights $\gl_{i}$ and $\gl_{i}-\gc$ respectively.
The observations of the previous Subsection give the following.
\bl\label{sea} The leading term of the coefficient of $e_{-\gc}v_\gl$ in $$e_{\ga_{1}}\ldots e_{\ga_r}e_{-\gc}e_{-\ga_{r}}\ldots e_{-\ga_1}v_\gl$$ is, up to a constant multiple, equal to $\prod_{i=1}^r (\gl,\ga_i).$ \el

\bl \label{alp}Suppose $\gc = w\gb$ with $w \in W_{\rm nonisotropic}$ and $\gb$ a simple isotropic root, and that $\gc$ is a simple root of some Borel. Then $q(w,\ga) =1$ for all $\ga \in N(w^{-1}),$ and
\[\{\gc-\ga_i| i=1, \ldots, r, (\gc,\ga_i^\vee) \neq 0\} = N(w^{-1}).
\]\el
\bpf The first statement can be checked on a case-by-case basis. The second is clearly true for $\gc$ simple.  Otherwise we can find a simple root $\ga$ such that $(\gc,\ga^\vee) >0$.  Then write $w= s_\ga u$ with $\ell(w) = \ell(u)+1$ and $u\in W_{\rm nonisotropic}$, and set $\gc' = s_\ga \gc$. The result then follows by induction on $\ell(w)$ and Equation (\ref{Nw}).  \epf
\noi
 There exists a unique (modulo a suitable left ideal in $U(\fg)$) \v Sapovalov element $\gth_\gc^{(i)}$ for the Borel subalgebra $\fb^{(i)}$ and polynomials $g_i(\gl), h_i(\gl)$ such that
\be \label{tlc}
e_{-\ga_i}\gth_\gc^{(i-1)}
e_{\ga_i}v_{i} = g_i(\gl)\gth_\gc^{(i)}v_{i}
\ee and
\be \label{tlc1}
e_{\ga_i}\gth_\gc^{(i)}e_{-\ga_i}v_{i-1} = h_i(\gl)\gth_\gc^{(i-1)}v_{i-1}
\ee
\bl \label{112} We have
\be g_i(\gl)h_i(\gl)=
(\gl +\gr,\ga_{i})(\gl+\gr-\gc,\ga_{i}).\ee
\el
\bpf Since $\gl_{i} =\gl_{i-1}-\ga_{i}$,
\begin{eqnarray}
\label{fff}
(\ga_{i},\gl_{i-1})(\ga_{i},\gl_{i-1}-\gc)\gth^{(i-1)}_\gc v_{i-1}
&=&e_{-\ga_{i}}e_{\ga_{i}} \gth_\gc^{(i)} e_{-\ga_{i}}e_{\ga_{i}}v_{i-1} .
\nonumber \\
&=&
h_i(\gl)e_{-\ga_{i}} \gth_\gc^{(i-1)} e_{\ga_{i}} v_{i-1}
\nonumber \\
&=&
g_i(\gl)h_i(\gl)\gth_\gc^{(i-1)}v_{i-1}
.\nonumber \end{eqnarray}
Now  each change of Borels is typical and  $\ga_{i}$ is simple isotropic for $\fb^{(i-1)}$, so
we have  $(\ga_{i-1},\gl_{i})$=$(\gl +\gr,\ga_{i}).$ The result follows.
\epf
 \bt \label{rot} Set $F(\gc) = \{i|1 \leq i \leq r \mbox{ and } (\gc,\ga_i)=0  \}$. There is a nonzero $c\in \ttk$ such that for all $\gl \in \cH_\gc$,

\[e_{\ga_{1}}\ldots e_{\ga_r}e_{-\gc}e_{-\ga_{r}}\ldots e_{-\ga_1}v_\gl= c\prod_{i \in F(\gc)}(\gl+\gr,\ga_i)
\gth_\gc v_\gl.\]
 \et
\bpf Since $v_i= e_{-\ga_i}v_{i-1},$ and $u_{i-1}=
e_{\ga_{i}}u_{i},$ Equation (\ref{tlc1}) and induction yield
  \be \label{pt} u_{i-1}=\prod_{j=i}^r h_j(\gl)\gth_\gc^{(i-1)}v_{i-1}. \ee Hence
  \be \label{nat} u_0 = \prod_{j=1}^r h_j(\gl) \gth_\gc v_{0}. \ee
  On the other hand, by Theorem \ref{aShap}  and Lemma \ref{alp} the leading term of the coefficient of $e_{-\gc}$ in $\gth_\gc v_\gl $ is given by
$\prod^r_{i=1, (\gc,\ga_i) \neq 0} (\gl,\ga_i)$
  up to a scalar multiple.
Therefore by comparing the coefficient of $e_{-\gc} v_\gl$ on both sides of (\ref{nat}), and using Lemma \ref{sea}, we have modulo terms of lower degree, that
\be \label{tat} \prod_{j=1}^r h_j(\gl)  = \prod_{i\in F(\gc)} (\gl+\gr,\ga_i). \ee
We note that the functions $\gl \lra (\gl+\gr,\ga_i)$ for $1 \leq i \leq r $ are linearly independent on
$\cH_\gc.$
It follows that  $h_j(\gl)$ is constant if
    $j \notin F(\gc),$ and $h_j(\gl) = (\gl+\gr,\ga_j)$
    if $j \in F(\gc).$
However we know from Lemma \ref{112} that if
$(\gl,\ga_j) = 0$ then  $h_j(\gl)$ divides  $(\ga_{j},\gl+\gr).$ Thus the result follows.
\epf \noi
Next suppose that $\fb'', \fb'$ are adjacent Borel subalgebras as in Equation {\rm (\ref{fggb})}, and that $d(\fb,\fb'')= d(\fb,\fb')+1$.  We can find a sequence of Borel subalgebras as in (\ref{distm}) such that $\fb'=\mathfrak{b}^{(i-1)}$ and
  $\fb''=\mathfrak{b}^{(i)}$ .  Adopting the notation of Equations (\ref{tlc}) and (\ref{tlc1}),
  we can now clarify the relationship between the \v Sapovalov elements
$\gth_\gc^{(i-1)}$ and $\gth_\gc^{(i)}$.
\bc With the above notation, we have up to constant multiples,
 \bi \itema If $(\gc,\ga_i) = 0$, then $h_i(\gl) =  g_i(\gl) = (\gl +\gr,\ga_{i})$. \itemb
If $(\gc,\ga_i) \neq 0$, then $h_i(\gl) = 1$ and $g_i(\gl) = (\gl +\gr,\ga_{i})(\gl +\gr-\gc,\ga_{i}).$
\ei\ec
\bpf This was shown in the course of the proof.\epf

\pagebreak




\end{document}